\pdfoutput=1
\RequirePackage{ifpdf}
\ifpdf
\documentclass[pdftex]{sigma}
\else
\documentclass{sigma}
\fi

\numberwithin{equation}{section}
\newtheorem{Theorem}{Theorem}[section]
\newtheorem{Corollary}[Theorem]{Corollary}
\newtheorem{Lemma}[Theorem]{Lemma}
\newtheorem{Proposition}[Theorem]{Proposition}
{\theoremstyle{definition}
\newtheorem{Definition}[Theorem]{Definition}
\newtheorem{Example}[Theorem]{Example}
\newtheorem{Remark}[Theorem]{Remark}
}
\DeclareMathOperator{\pt}{pt}
\DeclareMathOperator{\Mor}{Mor}
\DeclareMathOperator{\Hom}{Hom}
\DeclareMathOperator{\Id}{Id}
\DeclareMathOperator{\Spec}{Spec}

\begin{document}

\allowdisplaybreaks

\renewcommand{\PaperNumber}{022}

\FirstPageHeading

\ShortArticleName{The Real $K$-Theory of Compact Lie Groups}

\ArticleName{The Real $\boldsymbol{K}$-Theory of Compact Lie Groups}

\Author{Chi-Kwong FOK}

\AuthorNameForHeading{C.-K.~Fok}

\Address{Department of Mathematics, Cornell University, Ithaca, NY 14853, USA}
\Email{\href{mailto:ckfok@math.cornell.edu}{ckfok@math.cornell.edu}}
\URLaddress{\url{http://www.math.cornell.edu/~ckfok/}}

\ArticleDates{Received August 22, 2013, in f\/inal form March 06, 2014; Published online March 11, 2014}

\Abstract{Let~$G$ be a~compact, connected, and simply-connected Lie group, equipped with a~Lie group involution $\sigma_G$
and viewed as a~$G$-space with the conjugation action.
In this paper, we present a~description of the ring structure of the (equivariant) $KR$-theory of~$(G, \sigma_G)$ by
drawing on previous results on the module structure of the $KR$-theory and the ring structure of the equivariant
$K$-theory.}

\Keywords{$KR$-theory; compact Lie groups; Real representations; Real equivariant for\-ma\-lity}

\Classification{19L47; 57T10}

\vspace{-2mm}

\section{Introduction}
The complex $K$-theory of compact connected Lie groups with torsion-free fundamental groups was worked out by Hodgkin in
the 60s (cf.~\cite{Ho} and Theorem~\ref{hodgkin}).
It asserts that the $K$-theory ring is the $\mathbb{Z}_2$-graded exterior algebra over $\mathbb{Z}$ on the module of
primitive elements, which are of degree $-1$ and associated with the representations of the Lie group.
For the elegant proof of Hodgkin's result in the special case where~$G$ is simply-connected, see~\cite{At2}.

In~\cite{At3}, Atiyah introduced $KR$-theory, which is basically a~version of topological $K$-theory for the category of
Real spaces, i.e.\ topological spaces equipped with an involution.
$KR$-theory can be regarded as a~hybrid of $KO$-theory, complex $K$-theory and $KSC$-theory (cf.~\cite{At3}).
One can also consider equivariant $KR$-theory, where a~certain compatibility condition between the group action and the
involution is assumed.
For def\/initions and some basic properties, see Def\/initions~\ref{krtheory} and~\ref{eqkrtheory}~\cite{At3, AS}.

Since Hodgkin's work, there have appeared two kinds of generalizations of $K$-theory of compact Lie groups.
The f\/irst such is $KR$-theory of compact Lie groups, which was f\/irst studied by Seymour (cf.~\cite{Se}).
He obtained the $KR^*(\pt)$-module structure of $KR^*(G)$, where~$G$ is a~compact, connected and simply-connected
Lie group equipped with a~Lie group involution, using his structure theorem of $KR$-theory
of a~certain type of spaces
(cf.\ Theorems~\ref{strthm},~\ref{krtheorymodofg} and~\ref{spkrtheorymodofg}).
He was unable to obtain a~complete description of the ring structure, however, and could only make some conjectures
about it.
In~\cite{Bo}, Bousf\/ield determined functorially the united 2-adic $K$-cohomology algebra of any compact, connected and
simply-connected Lie group, which includes the 2-adic $KO$-cohomology algebra, and hence extended Seymour's results in
the 2-adic case, if the Lie group involution is taken into account appropriately.

The second one is the equivariant $K$-theory of compact Lie groups.
In~\cite{BZ}, Brylinski and Zhang showed that, for a~compact connected Lie group~$G$ with torsion free fundamental
group, the equivariant $K$-theory, $K_G^*(G)$, where~$G$ acts on itself by conjugation, is isomorphic to the ring of
Grothendieck dif\/ferentials of the complex representation ring $R(G)$ (for def\/inition, see~\cite{BZ}).
It is noteworthy that~$G$ satisf\/ies the property of being weakly equivariantly formal \emph{\`a la} Harada and Landweber
(cf.\ Def\/inition~4.1 of~\cite{HL} and Remark~\ref{weakequivformality}).

In this paper, based on the previous works of Seymour's and Brylinski--Zhang's and putting both the Real and equivariant
structures together, we obtain a~description of the ring structure of the equivariant $KR$-theory of any compact,
connected and simply-connected Lie group, which is recorded in Theorem~\ref{mainthm2}.
We express the ring structure using relations of generators associated to Real representations of~$G$ of real, complex
and quaternionic type (with respect to the Lie group involution).
Our main contribution is twofold.
First, we observe that the conditions of Seymour's structure theorem are an appropriate candidate for def\/ining the
notion of `Real formality' in analogy to weakly equivariant formality.
These notions together prompt us to introduce the def\/inition of `Real equivariant formality', which leads to a~structure
theorem for equivariant $KR$-theory (Theorem~\ref{equivstrthm}).
Any compact, connected and simply-connected Real Lie group falls under the category of Real equivariantly formal spaces
and Theorem~\ref{equivstrthm} enables us to obtain a~preliminary description of the equivariant $KR$-theory as an
algebra over the coef\/f\/icient ring.
Second, inspired by Seymour's conjecture, we obtain the squares of the real and quaternionic type generators, which in
addition to other known relations complete the description of the ring structure.
These squares are non-zero 2-torsions in general.
Hence the equivariant $KR$-theory in general is not a~ring of Grothendieck dif\/ferentials, as in the case for equivariant
$K$-theory.
Despite this, we remark that, in certain cases, if we invert~2 in the equivariant $KR$-theory ring, then the result is
an exterior algebra over the localized coef\/f\/icient ring of equivariant $KR$-theory.

The organization of this paper is as follows.
In Section~\ref{Section2}, we review the (equivariant) $K$-theory of compact connected Lie groups with torsion-free fundamental
groups, Real representation rings $RR(G)$, $KR$-theory, and the main results in~\cite{Se}.
In Section~\ref{Section3} we give a~description of the coef\/f\/icient ring $KR^*_G(\pt)$.
Section~\ref{Section4} is devoted to proving Theorems~\ref{equivstrthm} and~\ref{mainthm2}, the main results of this paper, which
give a~full description of the ring structure of $KR^*_G(G)$.
In Section~\ref{Section5} we apply the main results to obtain the ring structure of the ordinary $KR$-theory of compact Lie groups,
thereby conf\/irming Seymour's conjecture on the squares of the real type generators and disproving that on the squares of
the quaternionic type generators.
We also work out several examples, one of which shows how equivariant $KR$-theory tells apart (while $K$-theory cannot)
the ring structures of two equivariant $KR$-theory rings of~$G$, where in one case~$G$ acts on itself trivially, while
in another~$G$ acts by conjugation.

Throughout this paper, we follow the convention in~\cite{At3} and~\cite{AS} in reference to $KR$-theory.
In particular, we use the word `Real' (with a~capital R) in all contexts involving involutions, so as to avoid confusion
with the word `real' with the usual meaning.
For example, `Real $K$-theory' is used interchangeably with $KR$-theory, whereas `real $K$-theory' means $KO$-theory.

\vspace{-1mm}

\section{Preliminaries}\label{Section2}

Throughout this section,~$G$ denotes any compact Lie groups, and $X, Y, \dots$ any f\/inite $CW$-complexes unless
otherwise specif\/ied.

\vspace{-1mm}

\subsection[(Equivariant) $K$-theory of compact connected Lie groups]{(Equivariant)
$\boldsymbol{K}$-theory of compact connected Lie groups}

Recall that the functor $K^{-1}$ is represented by $U(\infty):=\lim\limits_{n\to\infty}U(n)$, i.e.\
for any $X$, $K^{-1}(X)$ is the abelian group of homotopy classes of maps $[X, U(\infty)]$.
Such a~description of $K^{-1}$ leads to the following
\begin{Definition}
\label{betamap}
Let $\delta: R(G)\to K^{-1}(G)$ be the group homomorphism which sends any complex~$G$-representation $\rho$ to the
homotopy class of $i\circ \rho$, where $\rho$ is regarded as a~continuous map from~$G$ to $U(n)$ and $i:
U(n)\hookrightarrow U(\infty)$ is the standard inclusion.
\end{Definition}

{\samepage In fact, any element in $K^{-q}(X)$ can be represented by a~complex of vector bundles on $X\times \mathbb{R}^q$, which
is exact outside $X\times\{0\}$ (cf.~\cite{At}).
This gives another interpretation of $\delta(\rho)$, as shown in the following proposition,
which we f\/ind useful for our
exposition.}

\begin{Proposition}
\label{-1_class}
If~$V$ is the underlying complex vector space of $\rho$, then $\delta(\rho)$ is represented by
\begin{gather*}
0\to G\times\mathbb{R}\times V \to G\times\mathbb{R}\times V\to 0,
\\
(g, t, v) \ \mapsto\  (g, t, -t\rho(g)v) \qquad\text{if} \ \ t\geq 0,
\\
(g, t, v) \ \mapsto\  (g, t, tv) \qquad\text{if}\ \  t\leq 0.
\end{gather*}
\end{Proposition}

\begin{Proposition}[\cite{Ho}]\label{noneqderivation}\quad
\begin{enumerate}\itemsep=0pt
\item[$1.$]
$\delta$ is a~derivation of $R(G)$ taking values in $K^{-1}(G)$ regarded as an $R(G)$-module whose module structure is
given by the augmentation map.
In other words, $\delta$ is a~group homomorphism from $R(G)$ to $K^{-1}(G)$ satisfying
\begin{gather}
\label{leibniz}
\delta(\rho_1\otimes\rho_2)=\dim(\rho_1)\delta(\rho_2)+\dim(\rho_2)\delta(\rho_1).
\end{gather}
\item[$2.$]
If $I(G)$ is the augmentation ideal of $R(G)$, then $\delta(I(G)^2)=0$.
\end{enumerate}
\end{Proposition}

The main results of~\cite{Ho} are stated in the following
\begin{Theorem}
\label{hodgkin}
Let~$G$ be a~compact connected Lie group with torsion-free fundamental group.
Then
\begin{enumerate}\itemsep=0pt
\item[$1.$]
$K^*(G)$ is torsion-free.
\item[$2.$]
Let $J(G):=I(G)/I(G)^2$.
Then the map $\widetilde{\delta}: J(G)\to K^{-1}(G)$ induced by $\delta$ is well-defined, and
$K^*(G)=\bigwedge(\text{\rm Im}(\widetilde{\delta}))$.
\item[$3.$]
In particular, if~$G$ is compact, connected and simply-connected of rank $l$, then
$K^*(G)=\bigwedge_{\mathbb{Z}}(\delta(\rho_1), \dots, \delta(\rho_l))$, where $\rho_1, \dots, \rho_l$ are the
fundamental representations.
\end{enumerate}
\end{Theorem}

Viewing~$G$ as a~$G$-space via the adjoint action, one may consider the equivariant $K$-theo\-ry~$K_G^*(G)$.
Let $\Omega^*_{R(G)/\mathbb{Z}}$ be the ring of Grothendieck dif\/ferentials of $R(G)$ over $\mathbb{Z}$, i.e.\
the exterior algebra over $R(G)$ of the module of K\"ahler dif\/ferentials of $R(G)$ over $\mathbb{Z}$ (cf.~\cite{BZ}).
\begin{Definition}
\label{defdeltag}
Let $\varphi: \Omega_{R(G)/\mathbb{Z}}^*\to K_G^*(G)$ be the $R(G)$-algebra homomorphism def\/ined by the following
\begin{enumerate}\itemsep=0pt
\item[1)]
$\varphi(\rho_V):=[G\times V]\in K_G^*(G)$, where~$G$ acts on $G\times V$ by $g_0\cdot(g_1, v)=(g_0g_1g_0^{-1},
\rho_V(g_0)v)$,
\item[2)]
$\varphi(d\rho_V)\in K_G^{-1}(G)$ is the complex of vector bundles in Def\/inition~\ref{-1_class} where~$G$ acts on
$G\times \mathbb{R}\times V$ by $g_0\cdot(g_1, t, v)=(g_0g_1g_0^{-1}, t, \rho_V(g_0)v)$,
\end{enumerate}
We also def\/ine $\delta_G: R(G)\to K^{-1}_G (G)$ by $\delta_G(\rho_V):=\varphi(d\rho_V)$.
\end{Definition}
\begin{Remark}
The def\/inition of $\delta_G(\rho_V)$ given in~\cite{BZ}, where the middle map of the complex of vector bundles is
def\/ined to be $(g, t, v)\mapsto(g, t, t\rho_V(g)v)$ for all $t\in\mathbb{R}$, is incorrect, as $\delta_G(\rho_V)$ so
def\/ined is actually 0 in $K_G^{-1}(G)$.
The def\/inition given in Def\/inition~\ref{defdeltag} is the correction made by Brylinski and relayed to us by one of the
referees.
\end{Remark}
\begin{Theorem}[\cite{BZ}]\label{eqderivation}\quad
\begin{enumerate}\itemsep=0pt
\item[$1.$]
$\delta_G$ is a~derivation of $R(G)$ taking values in the $R(G)$-module $K^{-1}_G(G)$, i.e.\
$\delta_G$ satisfies
\begin{gather}
\label{eqleibniz}
\delta_G(\rho_1\otimes\rho_2)=\rho_1\cdot\delta_G(\rho_2)+\rho_2\cdot\delta_G(\rho_1).
\end{gather}
\item[$2.$]
Let~$G$ be a~compact connected Lie group with torsion-free fundamental group.
Then $\varphi$ is an $R(G)$-algebra isomorphism.
\end{enumerate}
\end{Theorem}

\begin{Remark}\label{weakequivformality}\quad
\begin{enumerate}\itemsep=0pt
\item
Although the def\/inition of $\delta_G(\rho_V)$ given by Brylinski--Zhang in~\cite{BZ} is incorrect, their proof of
Theorem~\ref{eqderivation} can be easily corrected by using the correct def\/inition as in Def\/inition~\ref{defdeltag},
which does not af\/fect the validity of the rest of their arguments, and Theorem~\ref{eqderivation} still stands.
\item
In~\cite{HL}, a~$G$-space $X$ is def\/ined to be \emph{weakly equivariantly formal} if the map
$K^*_G(X)\otimes_{R(G)}\mathbb{Z}\to K^*(X)$ induced by the forgetful map is a~ring isomorphism, where $\mathbb{Z}$ is
viewed as an $R(G)$-module through the augmentation homomorphism.
Theorem~\ref{hodgkin} and~\ref{eqderivation} imply that~$G$ is weakly equivariantly formal if it is connected with
torsion-free fundamental group.
We will make use of this property in computing the equivariant $KR$-theory of~$G$.
\item
Let $f: K_G^*(G)\to K^*(G)$ be the forgetful map.
Note that $f(\varphi(\rho))=\dim(\rho)$ and $f(\varphi(d\rho))=\delta(\rho)$.
Applying $f$ to equation~\eqref{eqleibniz} in Theorem~\ref{eqderivation}, we get equation~\eqref{leibniz} in
Proposition~\ref{noneqderivation}.
\item
$I(G)$, being a~prime ideal in $R(G)$, can be thought of as an element in $\Spec R(G)$, and
$K^*(G)\cong\bigwedge_{\mathbb{Z}}T^*_{I(G)}\Spec R(G)$, $K_G^*(G)\cong\bigwedge_{R(G)}T^*_{I(G)}\Spec R(G)$.
\end{enumerate}
\end{Remark}
\subsection{Real representation rings}
This section is an elaboration of the part on Real representation rings in~\cite{AS} and~\cite{Se}.
Since the results in this subsection can be readily generalized from those results concerning the special case
where $\sigma_G$ is trivial, we omit most of the proofs and refer the reader to any standard text on representation theory of
Lie groups, e.g.~\cite{B}
and~\cite{BtD}.
\begin{Definition}
A Real Lie group is a~pair $(G,\sigma_G)$ where~$G$ is a~Lie group and $\sigma_G$ a~Lie group involution on it.
A \emph{Real representation}
$V$ of a~Real Lie group $(G,\sigma_G)$ is a~f\/inite-dimensional complex representation of~$G$ equipped with an anti-linear involution~$\sigma_V$ such
that $\sigma_V(gv)=\sigma_G(g)\sigma_V(v)$.
Let $\mathcal{R}\text{ep}_\mathbb{R}(G,
\sigma_G)$ be the category of Real representations of $(G, \sigma_G)$.
A morphism between~$V$ and $W\in \mathcal{R}
ep_\mathbb{R}(G,
\sigma_G)$ is a~linear transformation from~$V$ to $W$ which commute with~$G$ and respect~$\sigma_V$ and $\sigma_W$.
We denote $\Mor
(V, W)$ by $\Hom_G(V, W)^{(\sigma_V, \sigma_W)}$.
An irreducible Real representation is an irreducible object in $\mathcal{R}\text{ep}_\mathbb{R}(G,
\sigma_G)$.
The Real representation ring of $(G, \sigma_G)$, denoted by $RR(G, \sigma_G)$, is the Grothendieck group of $\mathcal{R}
ep_\mathbb{R}(G,
\sigma_G)$, with multiplication being tensor product over $\mathbb{C}$.
Sometimes we will omit the notation $\sigma_G$ if there is no danger of confusion about the Lie group involution.
\end{Definition}
\begin{Remark}
Let~$V$ be an irreducible Real representation of~$G$.
Then $\Hom_G(V, V)^{\sigma_V}$ must be a~real associative division algebra which, according to Frobenius' theorem, is isomorphic to either
$\mathbb{R}$, $\mathbb{C}$ or $\mathbb{H}$.
Following the convention of~\cite{AS}, we call $\Hom_G(V, V)^{\sigma_V}$ the \emph{commuting field} of~$V$.
\end{Remark}

\begin{Definition}
If~$V$ is an irreducible Real representation of~$G$, then we say~$V$ is of \emph{real, complex or quaternionic type}
according as the commuting f\/ield is isomorphic to $\mathbb{R}$, $\mathbb{C}$ or $\mathbb{H}$.
Let $RR(G, \mathbb{F})$ be the abelian group generated by the isomorphism classes of irreducible Real representations
with $\mathbb{F}$ as the commuting f\/ield.
\end{Definition}
\begin{Remark}
$RR(G)\cong RR(G, \mathbb{R})\oplus RR(G, \mathbb{C})\oplus RR(G, \mathbb{H})$ as abelian groups.
\end{Remark}
\begin{Definition}
Let~$V$ be a~$G$-representation.
We use $\sigma_G^*V$ to denote the~$G$-representation with the same underlying vector space where the~$G$-action is twisted by
$\sigma_G$, i.e.\
$\rho_{\sigma_G^*V}
(g)v=\rho_V(\sigma_G(g))v$.
We will use $\overline{\sigma_G^*}$ to denote the map on $R(G)$ def\/ined by $[V]\mapsto [\sigma_G^*\overline{V}]$.
\end{Definition}

\begin{Proposition}
If~$V$ is a~complex~$G$-representation, and there exists $f\in\Hom_G(V,\sigma_G^*\overline{V})$
such that $f^2=\Id_V$, then~$V$ is a~Real representation of~$G$ with $f$ as the anti-linear involution~$\sigma_V$.
\end{Proposition}
\begin{proof}
If $f\in\Hom_G(V,
\sigma_G^*\overline{V})$, then it is anti-linear on~$V$ and $f(gv)=
\sigma_G(g)f(v)$ for $g\in G$ and $v\in V$.
The assumption that $f^2=\Id
_V$ just says that $f$ is an involution.
So~$V$ together with $\sigma_V=f$ is a~Real representation of~$G$.
\end{proof}
\begin{Proposition}
Let~$V$ be an irreducible complex representation of~$G$ and suppose that $V\cong
\sigma_G^*\overline{V}$.
Let $f\in\Hom_G(V,
\sigma_G^*\overline{V})$.
Then
\begin{enumerate}\itemsep=0pt
\item[$1.$]
$f^2=k\Id_V$ for some $k\in\mathbb{R}$.
\item[$2.$]
There exists $g\in\Hom_G(V,
\sigma_G^*\overline{V})$ such that $g^2=\Id_V$ or $g^2=-\Id_V$.
\end{enumerate}
\end{Proposition}
\begin{proof}
Note that $f^2\in\Hom_G(V, V)$.
By Schur's lemma, $f^2=k\Id_V$ for some $k\in\mathbb{C}$.
On the other hand,
\begin{gather*}
f(kv)=f(f(f(v)))=kf(v).
\end{gather*}
But $f$ is an anti-linear map on~$V$.
It follows that $k=\overline{k}$ and hence $k\in\mathbb{R}$.
For part (2), we may f\/irst simply pick an isomorphism $f\in\Hom_G(V,
\sigma_G^*\overline{V})$.
Then $f^2=k\Id_V$ for some $k\in\mathbb{R}^\times$.
Schur's lemma implies that any $g\in\Hom_G(V,
\sigma_G^*\overline{V})$ must be of the form $g=cf$ for some $c\in\mathbb{C}$.
Then $g^2=cf\circ cf=c\overline{c}f^2=|c|^2k\Id_V$.
Consequently, if $k$ is positive, we choose $c=\frac{1}{\sqrt{k}}$ so that $g^2=\Id_V$; if $k$ is negative, we choose $c=\frac{1}{\sqrt{-k}}$ so that $g^2=-\Id_V$.
\end{proof}

\begin{Proposition}
\label{realrepclassification}
Let~$V$ be an irreducible Real representation of~$G$.
\begin{enumerate}\itemsep=0pt
\item[$1.$]
The commuting field of~$V$ is isomorphic to $\mathbb{R}$ iff~$V$ is an irreducible complex representation and there
exists $f\in\Hom_G(V,
\sigma_G^*\overline{V})$ such that $f^2=\Id_V$.
\item[$2.$]
The commuting field of~$V$ is isomorphic to $\mathbb{C}$ iff $V\cong W\oplus
\sigma_G^*\overline{W}$ as complex~$G$-representa\-tions,
where $W$ is an irreducible complex~$G$-representation and $W\ncong
\sigma_G^*\overline{W}$, and $\sigma_V(w_1, w_2)$ $=(w_2, w_1)$.
\item[$3.$]
The commuting field of~$V$ is isomorphic to $\mathbb{H}$ iff $V\cong W\oplus
\sigma_G^*\overline{W}$ as complex~$G$-repre\-sen\-ta\-tions,
where $W$ is an irreducible complex~$G$-representation and there exists
$f\in\Hom_G(V,\sigma_G^*\overline{V})$ such that $f^2=-\Id_V$, and $\sigma_V(w_1, w_2)=(w_2, w_1)$.
\end{enumerate}
\end{Proposition}

\begin{proof}
One can easily establish the above proposition
by modifying the proof of Proposition~3 in Appendix~2 of~\cite{B}, which
is a~special case of the above proposition
where $\sigma_G$ is trivial.
\end{proof}
\begin{Proposition}\label{realforgetfulinj}\quad
\begin{enumerate}\itemsep=0pt
\item[$1.$]
The map
$i: RR(G)\to R(G)$ which forgets the Real structure is injective.
\item[$2.$]
Any complex~$G$-representation~$V$ which is a~Real representation can only possess a~unique Real structure up to
isomorphisms of Real~$G$-representations.
\end{enumerate}
\end{Proposition}
\begin{proof}
Let $\rho: R(G)\to RR(G)$ be the map
\begin{gather*}
[V]\mapsto \big[\big(V\oplus\iota_G^*\overline{V},\sigma_{V\oplus\iota_G^*\overline{V}}\big)\big],
\end{gather*}
where $\sigma_{V\oplus\iota_G^*\overline{V}}(u, w)=(w, u)$.
Let $[(V,
\sigma_V)]\in RR(G)$.
Then
\begin{gather*}
\rho\circ i([(V,\sigma_V)])=\big[\big(V\oplus\iota_G^*\overline{V},\sigma_{V\oplus\iota_G^*\overline{V}}\big)\big].
\end{gather*}
We claim that $[(V\oplus V,\sigma_V\oplus\sigma_V)]=[(V\oplus\iota_G^*\overline{V},
\sigma_{V\oplus\iota_G^*\overline{V}})]$, which is easily seen to be true because of the Real~$G$-representation isomorphism
\begin{gather*}
\begin{split}
& f: \ V\oplus \iota_G^*\overline{V}  \to  V\oplus V,
\\
& \phantom{f:{}} \ (u, w)  \mapsto  (u+\sigma_V(w), i(u-\sigma_V(w))).
\end{split}
\end{gather*}
It follows that $\rho\circ i$ amounts to multiplication by 2 on $RR(G)$, and is therefore injective because $RR(G)$ is
a~free abelian group generated by irreducible Real representations.
Hence $i$ is injective.
(2) is simply a~restatement of (1).
\end{proof}
Proposition~\ref{realforgetfulinj} makes it legitimate to regard $RR(G)$ as a~subring of $R(G)$.
From now on we view $[V]\in R(G)$ as an element in $RR(G)$ if~$V$ possesses a~compatible Real structure.
\begin{Proposition}
\label{formclassification}
Let~$G$ be a~compact Real Lie group.
Let~$V$ be an irreducible complex representation of~$G$.
Then
\begin{enumerate}\itemsep=0pt
\item[$1.$]
$[V]\in RR(G, \mathbb{R})$ iff there exists a~$G$-invariant symmetric nondegenerate bilinear form $B: V\times
\sigma_G^*V\to\mathbb{C}$.
\item[$2.$]
$[V\oplus V]\in RR(G, \mathbb{H})$ iff there exists a~$G$-invariant skew-symmetric nondegenerate bilinear form $B:
V\times
\sigma_G^*V\to\mathbb{C}$.
\item[$3.$]
$[V\oplus
\sigma_G^*\overline{V}
]\in RR(G, \mathbb{C})$ iff there does not exist any~$G$-invariant nondegenerate bilinear form on $V\times
\sigma_G^*V$.
\end{enumerate}
\end{Proposition}
\begin{proof}
By Proposition~\ref{realrepclassification},~$V$ is a~Real representation of real type if\/f there exists
$f\in\Hom(V,
\sigma_G^*\overline{V})^G$ such that $f^2=\Id_V$.
One can def\/ine a~bilinear form $B: V\times
\sigma_G^*V\to\mathbb{C}$ by
\begin{gather}
\label{bilinearformtype}
B(v_1, v_2)=\langle v_1, f(v_2)\rangle,
\end{gather}
where $\langle \,, \,\rangle$ is a~$G$-invariant Hermitian inner product on~$V$.
It can be easily seen that $B$ is~$G$-invariant, symmetric and non-degenerate.
Conversely, given a~$G$-invariant symmetric non-degenerate bilinear form on $V\times
\sigma_G^*V$ and using equation~\eqref{bilinearformtype}, we can def\/ine $f\in\Hom(V,
\sigma_G^*\overline{V})^G$, which squares to identity.
Part (2) and (3) follow similarly.
\end{proof}
Proposition~\ref{formclassification} leads to the following

\begin{Definition}
\label{cplxrepclassification}
Let~$G$ be a~compact Real Lie group.
Let~$V$ be an irreducible complex representation of~$G$.
Def\/ine, with respect to $\sigma_G$,
\begin{enumerate}\itemsep=0pt
\item
$V$ to be of real type if there exists a~$G$-invariant symmetric nondegenerate bilinear form $B: V\times
\sigma_G^*V\to\mathbb{C}$.
\item
$V$ to be of quaternionic type if there exists a~$G$-invariant skew-symmetric nondegenerate bilinear form $B: V\times
\sigma_G^*V\to\mathbb{C}$.
\item
$V$ to be of complex type if $V\ncong
\sigma_G^*\overline{V}$.
\end{enumerate}
The abelian group generated by classes of irreducible complex representation of type $\mathbb{F}$ is denoted by $R(G,
\mathbb{F})$.
\end{Definition}

\begin{Definition}
If~$V$ is a~complex~$G$-representation equipped with an anti-linear endomor\-phism~$J_V$ such that $J_V(gv)=
\sigma_G(g)J(v)$ and $J^2=-\Id
_V$, then we say~$V$ is a~\emph{Quaternionic representation} of~$G$.
Let $\mathcal{R}\text{ep}_\mathbb{H}(G)$ be the category of Quaternionic representations of~$G$.
A~morphism between~$V$ and $W\in \mathcal{R}\text{ep}_\mathbb{H}(G)$ is a~linear transformation from~$V$ to $W$ which commutes
with~$G$ and respect $J_V$ and $J_W$.
We denote $\Mor(V, W)$ by $\Hom_G(V, W)^{(J_V, J_W)}$.
An irreducible Quaternionic representation is an irreducible object in $\mathcal{R}\text{ep}_\mathbb{H}(G,
\sigma_G)$.
The Quaternionic representation group of~$G$, denoted by $RH(G)$, is the Grothendieck group of $\mathcal{R}\text{ep}_\mathbb{H}(G)$.
Let $RH(G, \mathbb{F})$ be the abelian group generated by the isomorphism classes of irreducible Quaternionic
representations with $\mathbb{F}$ as the commuting f\/ield.
\end{Definition}
\begin{Remark}
The tensor product of two Quaternionic representations~$V$ and $W$ is a~Real representation as $J_V\otimes J_W$ is an
anti-linear involution which is compatible with $\sigma_G$.
Similarly the tensor product of a~Real representation and a~Quaternionic representation is a~Quaternionic
representation.
To put it succinctly, $RR(G)\oplus RH(G)$ is a~$\mathbb{Z}
_2$-graded ring, with $RR(G)$ being the degree $0$ piece and $RH(G)$ the degree $-1$ piece.
Later on we will assign $RH(G)$ with a~dif\/ferent degree so as to be compatible with the description of the coef\/f\/icient
ring $KR^*_G(\pt)$.
\end{Remark}
\begin{Proposition}
\label{quatrepclassification}
$RH(G)$, as an abelian group, is generated by the following
\begin{enumerate}\itemsep=0pt
\item[$1.$]
$[V\oplus
\sigma_G^*\overline{V}]$ where~$V$ is an irreducible complex representation of real type with $J(u, w)=(-w, u)$.
Its commuting field is $\mathbb{H}$.
\item[$2.$]
$[V\oplus\sigma_G^*\overline{V}]$, where~$V$ is an irreducible complex representation of complex type with $J(u, w)=(-w, u)$.
Its commuting field is $\mathbb{C}$.
\item[$3.$]
$[V]$, where~$V$ is an irreducible complex representation of quaternionic type.
Its commuting field is $\mathbb{R}$.
\end{enumerate}
\end{Proposition}
\begin{proof}
The proof proceeds in the same fashion as does the proof for Proposition~\ref{realrepclassification}.
\end{proof}
\begin{Corollary}\label{3by3table}\quad
\begin{enumerate}\itemsep=0pt
\item[$1.$]
We have that $RR(G, \mathbb{R})\cong RH(G, \mathbb{H})\cong R(G, \mathbb{R})$, $RH(G, \mathbb{R})\cong RR(G,\mathbb{H})$
$\cong R(G, \mathbb{H})$ and $RR(G, \mathbb{C})\cong RH(G, \mathbb{C})$ as abelian groups.
\item[$2.$]
$RR(G)$ is isomorphic to $RH(G)$ as abelian groups.
\end{enumerate}
\end{Corollary}

\begin{proof}
The result follows easily from Proposition~\ref{realrepclassification}, Def\/inition~\ref{cplxrepclassification} and
Proposition~\ref{quatrepclassification}.
\end{proof}
\begin{Proposition}\label{quatforgetfulinj}\quad
\begin{enumerate}\itemsep=0pt
\item[$1.$]
The map $j: RH(G)\to R(G)$ which forgets the Quaternionic structure is injective.
\item[$2.$]
Any complex~$G$-representation which is a~Quaternionic representation can only possess a~unique Quaternionic structure
up to isomorphisms of Quaternionic~$G$-representations.
\end{enumerate}
\end{Proposition}
\begin{proof}
The proof proceeds in the same fashion as does the proof for Proposition~\ref{realforgetfulinj}.
It suf\/f\/ices to show that, if $\eta: R(G)\to RH(G)$ is the map $[V]\mapsto [(V\oplus
\sigma_G^*\overline{V},
\sigma_{V\oplus\iota_G^*\overline{V}})]$, then $\eta\circ j$ amounts to multiplication by 2, i.e.\
$(V\oplus\sigma_G^*\overline{V},\sigma_{V\oplus\sigma_G^*\overline{V}})\cong(V\oplus V, J\oplus J)$, where $\sigma_{V\oplus\sigma_G^*\overline{V}}(u, w)=(-w, u)$.
This is true because of the Quaternionic~$G$-representation isomorphism
\begin{gather*}
f: \  V\oplus\iota_G^*\overline{V}  \to   V\oplus V,\\
\hphantom{f: {}}  \ (u, w) \mapsto  (u+Jw, i(u-Jw)). \tag*{\qed}
\end{gather*}
\renewcommand{\qed}{}
\end{proof}

\begin{Example}
We shall illustrate the similarities and dif\/ferences of the various aforementioned representation groups with an
example.
Let $G=Q_8\times C_3$, the direct product of the quaternion group and the cyclic group of order 3, equipped with the
trivial involution.
There are 5 irreducible complex representations of $Q_8$, namely, the 4 1-dimensional representations which become
trivial on restriction to the center $Z$ of $Q_8$ and descend to the 4 1-dimensional representations of $Q_8/Z\cong
\mathbb{Z}_2\oplus\mathbb{Z}_2$, and the 2-dimensional faithful representation.
We denote these representations by $1_{Q_8}$, $\rho_{(0, 1)}$, $\rho_{(1, 0)}$, $\rho_{(1, 1)}$ and~$\rho_Q$ respectively.
Similarly, we let~$1_{C_3}$,~$\rho_\zeta$ and~$\rho_{\zeta^2}$ be the three 1-dimensional complex irreducible
representations of~$C_3$.
It can be easily seen that
\begin{gather*}
R(Q_8, \mathbb{R})=\mathbb{Z}\cdot[1_{Q_8}]\oplus\mathbb{Z}\cdot[\rho_{(0, 1)}]\oplus\mathbb{Z}\cdot[\rho_{(1,
0)}]\oplus\mathbb{Z}\cdot[\rho_{(1, 1)}],
\\
R(Q_8, \mathbb{C})=0,
\\
R(Q_8, \mathbb{H})=\mathbb{Z}\cdot[\rho_Q],
\\
R(C_3, \mathbb{R})=\mathbb{Z}\cdot[1_{C_3}],
\\
R(C_3, \mathbb{C})=\mathbb{Z}\cdot[\rho_\zeta]\oplus\mathbb{Z}\cdot[\rho_{\zeta^2}],
\\
R(C_3, \mathbb{H})=0.
\end{gather*}
It follows that
\begin{gather*}
R(G, \mathbb{R})=\bigoplus_{x\in \{1_{Q_8}, \rho_{(0, 1)}, \rho_{(1, 0)}, \rho_{(1,
1)}\}}\mathbb{Z}\cdot[x\widehat{\otimes}1_{C_3}]\cong\mathbb{Z}^4,
\\
R(G, \mathbb{C})=\bigoplus_{
\begin{subarray}
{c}{x\in\{1_{Q_8}, \rho_{(0, 1)}, \rho_{(1, 0)}, \rho_{(1, 1)}, \rho_Q\}}
\\
{y\in\{\rho_\zeta, \rho_{\zeta^2}\}}
\end{subarray}}\mathbb{Z}\cdot[x\widehat{\otimes}y]\cong\mathbb{Z}^{10},
\\
R(G, \mathbb{H})=\mathbb{Z}\cdot[\rho_Q\widehat{\otimes}1_{C_3}]\cong\mathbb{Z},
\\
RR(G, \mathbb{R})=\bigoplus_{x\in \{1_{Q_8}, \rho_{(0, 1)}, \rho_{(1, 0)}, \rho_{(1,
1)}\}}\mathbb{Z}\cdot[x\widehat{\otimes}1_{C_3}]\cong\mathbb{Z}^4,
\\
RR(G, \mathbb{C})=\bigoplus_{x\in \{1_{Q_8}, \rho_{(0, 1)}, \rho_{(1, 0)}, \rho_{(1, 1)},
\rho_Q\}}\mathbb{Z}\cdot[x\widehat{\otimes}\rho_\zeta\oplus\overline{x}\widehat{\otimes}\rho_{\zeta^2}]\cong\mathbb{Z}^5,
\\
RR(G, \mathbb{H})=\mathbb{Z}\cdot[\rho_Q\widehat{\otimes}1_{C_3}\oplus
\overline{\rho_Q}\widehat{\otimes}1_{C_3}]\cong\mathbb{Z},
\\
RH(G, \mathbb{R})=\mathbb{Z}\cdot[\rho_Q\widehat{\otimes}1_{C_3}]\cong\mathbb{Z},
\\
RH(G, \mathbb{C})=\bigoplus_{x\in \{1_{Q_8}, \rho_{(0, 1)}, \rho_{(1, 0)}, \rho_{(1, 1)},
\rho_Q\}}\mathbb{Z}\cdot[x\widehat{\otimes}\rho_\zeta\oplus\overline{x}\widehat{\otimes}\rho_{\zeta^2}]\cong\mathbb{Z}^5,
\\
RH(G, \mathbb{H})=\bigoplus_{x\in \{1_{Q_8}, \rho_{(0, 1)}, \rho_{(1, 0)}, \rho_{(1,
1)}\}}\mathbb{Z}\cdot[x\widehat{\otimes}1_{C_3}\oplus\overline{x}\widehat{\otimes}1_{C_3}]\cong\mathbb{Z}^4.
\end{gather*}
Some representations above should be equipped with suitable Real or Quaternionic structures given in
Propositions~\ref{realrepclassification} and~\ref{quatrepclassification}.
For example, the Real structure of $\rho_Q\widehat{\otimes}1_{C_3}\oplus\overline{\rho_Q}\widehat{\otimes}1_{C_3}$ in
$RR(G, \mathbb{H})$ is given by swapping the two coordinates.
\end{Example}

\subsection[$KR$-theory]{$\boldsymbol{KR}$-theory}

$KR$-theory was f\/irst introduced by Atiyah in~\cite{At3} and used to derive the 8-periodicity of $KO$-theory from the
2-periodicity of complex $K$-theory.
$KR$-theory was motivated by the index theory of real elliptic operators.
\begin{Definition}\label{krtheory}\quad
\begin{enumerate}\itemsep=0pt
\item
A \emph{Real space} is a~pair $(X,
\sigma_X)$ where $X$ is a~topological space equipped with an involutive homeomorphism $\sigma_X$, i.e.\
$\sigma_X^2=\Id
_X$.
We will sometimes suppress the notation $\sigma_X$ and simply use $X$ to denote the Real space,
if there is no danger of confusion about the involutive
homeomorphism.
A \emph{Real pair}
is a~pair $(X, Y)$ where~$Y$ is a~closed subspace of~$X$ invariant under~$\sigma_X$.
\item
Let $\mathbb{R}^{p, q}$ be the Euclidean space $\mathbb{R}^{p+q}$
equipped with the involution which is identity on the f\/irst~$q$
coordinates and negation on the last $p$-coordinates.
Let $B^{p, q}$ and $S^{p, q}$ be the unit ball and sphere in $\mathbb{R}^{p, q}$ with the inherited involution.
\item
A~\emph{Real vector bundle}
(to be distinguished from the usual real vector bundle) over $X$ is a~complex vector bundle
$E$ over $X$ which itself is also a~Real space with involutive homeomorphism $\sigma_E$ satisfying
\begin{enumerate}\itemsep=0pt
\item
$\sigma_X\circ p=p\circ\sigma_E$, where $p: E\to X$ is the projection map,
\item
$\sigma_E$ maps $E_x$ to $E_{\sigma_X(x)}$ anti-linearly.
\end{enumerate}
A \emph{Quaternionic vector bundle}
(to be distinguished from the usual quaternionic vector bundle) over $X$ is a~complex
vector bundle $E$ over $X$ equipped with an anti-linear lift $\sigma_E$ of $\sigma_X$ such that $\sigma_E^2=-\Id
_E$.
\item
Let $X$ be a~Real space.
The ring $KR(X)$ is the Grothendieck group of the isomorphism classes of Real vector bundles over $X$, equipped with the
usual product structure induced by tensor product of vector bundles over $\mathbb{C}$.
The relative $KR$-theory for a~Real pair $KR(X, Y)$ can be similarly def\/ined.
In general, the graded $KR$-theory ring of the Real pair $(X, Y)$ is given by
\begin{gather*}
KR^*(X, Y):=\bigoplus_{q=0}^7 KR^{-q}(X, Y),
\end{gather*}
where
\begin{gather*}
KR^{-q}(X, Y):=KR\big(X\times B^{0,q}, X\times S^{0, q}\cup Y\times B^{0,q}\big).
\end{gather*}
The ring structure of $KR^*$ is extended from that of~$KR$, in a~way analogous to the case of complex $K$-theory.
The number of graded pieces, which is~8, is a~result of Bott periodicity for $KR$-theory (cf.~\cite{At3}).
\end{enumerate}
\end{Definition}

Note that when $\sigma_X=\Id_X$, then $KR(X)\cong KO(X)$.
On the other hand, if $X\times \mathbb{Z}_2$ is given the involution which swaps the two copies of $X$, then
$KR(X\times\mathbb{Z}_2)\cong K(X)$.
Also, if $X$ is equipped with the trivial involution, then $KR(X\times S^{2,0})\cong KSC(X)$, the Grothendieck group of
homotopy classes of self-conjugate bundles over $X$ (cf.~\cite{At3}).
In this way, it is natural to view $KR$-theory as a~unifying thread of $KO$-theory, $K$-theory and $KSC$-theory.

On top of the Real structure, we may further add compatible group actions and def\/ine equivariant $KR$-theory.
\begin{Definition}\label{eqkrtheory}\quad
\begin{enumerate}\itemsep=0pt
\item
A \emph{Real}~$G$-\emph{space} $X$ is a~quadruple $(X, G,
\sigma_X, \sigma_G)$ where a~group~$G$ acts on~$X$ and~$\sigma_G$ is an involutive automorphism of~$G$ such that
\begin{gather*}
\sigma_X(g\cdot x)=\sigma_G(g)\cdot\sigma_X(x).
\end{gather*}
\item
A \emph{Real}~$G$-\emph{vector bundle} $E$ over a~Real~$G$-space $X$ is a~Real vector bundle and a~$G$-bundle over $X$,
and it is also a~Real~$G$-space.
\item
In a~similar spirit, one can def\/ine equivariant $KR$-theory $KR_G^*(X, Y)$.
Notice that the~$G$-actions on $B^{0,q}$ and $S^{0,q}$ in the def\/inition of $KR_G^{-q}(X, Y)$ are trivial.
\end{enumerate}
\end{Definition}

\begin{Definition}\quad
\begin{enumerate}\itemsep=0pt
\item
Let $K^*(+)$ be the complex $K$-theory of a~point extended to a~$\mathbb{Z}_8$-graded algebra over
$K^0(\pt)\cong\mathbb{Z}$, i.e.\
$\displaystyle K^*(+)\cong\mathbb{Z}[\beta]\left/\beta^4-1\right.$.
Here $\beta\in K^{-2}(+)$ is the class of the reduced canonical bundle on $\mathbb{CP}^1\cong S^2$.
\item
Let $\overline{\sigma_X^*}$ be the map def\/ined on (equivariant) vector bundles on $X$ by $\overline{\sigma_X^*}
E:=
\sigma_X^*\overline{E}$.
The involution induced by $\overline{\sigma_X^*}$ on $K^*_G(X)$ is also denoted by $\overline{\sigma_X^*}$ for simplicity.
\end{enumerate}
\end{Definition}
In the following proposition, we collect, for reader's convenience, some basic results of $KR$-theory
(cf.~\cite[Section~2]{Se}), some of which are stated in the more general context of equivariant $KR$-theory.
\begin{Proposition}\label{krprelim}\quad
\begin{enumerate}\itemsep=0pt
\item[$1.$]
We have
\begin{gather*}
KR^*(\pt)\cong\mathbb{Z}[\eta, \mu]\big/\big(2\eta, \eta^3, \mu\eta, \mu^2-4\big) ,
\end{gather*}
where $\eta\in KR^{-1}(\pt)$, $\mu\in KR^{-4}(\pt)$ represents the reduced Hopf bundles of
$\mathbb{R}\mathbb{P}^1$ and $\mathbb{H}\mathbb{P}^1$ respectively.
\item[$2.$]
Let $c: KR^*_G(X)\to K^*_G(X)$ be the homomorphism which forgets the Real structure of Real vector bundles, and $r:
K^*_G(X)\to KR^*_G(X)$ be the realification map defined by $[E]\mapsto [E\oplus
\sigma_G^*\sigma_X^*\overline{E}
]$.
Then we have the following relations
\begin{enumerate}\itemsep=0pt
\item[$(a)$]
$c(1)=1$, $c(\eta)=0$, $c(\mu)=2\beta^2$, where $\beta\in K^{-2}(\pt)$ is the Bott class,
\item[$(b)$]
$r(1)=2$, $r(\beta)=\eta^2$, $r(\beta^2)=\mu$, $r(\beta^3)=0$,
\item[$(c)$]
$r(xc(y))=r(x)y$, $cr(x)=x+
\sigma_G^*\overline{\sigma_X^*}
x$ and $rc(y)=2y$ for $x\in K^*_G(X)$ and $y\in KR^*_G(X)$, where $K^*_G(X)$ is extended to a~$\mathbb{Z}_8$-graded
algebra by Bott periodicity.
\end{enumerate}
\end{enumerate}
\end{Proposition}
\begin{proof}
(1) is given in~\cite[Section~2]{Se}.
The proof of (2) is the same as in the nonequivariant case, which is given in~\cite{At3}.
\end{proof}

\begin{Definition}
A Quaternionic~$G$-vector bundle over a~Real space $X$ is a~complex vector bundle $E$ equipped with an anti-linear
vector bundle endomorphism $J$ on $E$ such that $J^2=-\Id_E$ and $J(g\cdot v)=
\sigma_G(g)\cdot J(v)$.
Let $KH^*_G(X)$ be the corresponding $K$-theory constructed using Quaternionic~$G$-bundles over $X$.
\end{Definition}
By generalizing the discussion preceding Lemma 5.2 in~\cite{Se} to the equivariant and graded setting, we def\/ine
a~natural transformation
\begin{gather*}
t: \ KH_G^{-q}(X)\to KR^{-q-4}_G(X)
\end{gather*}
which sends
\begin{gather*}
0\longrightarrow E_1
\stackrel{f}
{\longrightarrow} E_2\longrightarrow 0
\end{gather*}
to
\begin{gather*}
0\longrightarrow \pi^*(\mathbb{H}\otimes_\mathbb{C} E_1)\stackrel{g}
{\longrightarrow}\pi^*(\mathbb{H}\otimes_\mathbb{C}E_2)\longrightarrow 0,
\end{gather*}
where
\begin{enumerate}\itemsep=0pt
\item[1)]
$E_i$, $i=1, 2$ are equivariant Quaternionic vector bundles on $X\times\mathbb{R}^{0, q}$ equipped with the Quaternionic
structures $J_{E_i}$,
\item[2)]
$f$ is an equivariant Quaternionic vector bundle homomorphism which is an isomorphism outside $X\times\{0\}$,
\item[3)]
$\pi: X\times\mathbb{R}^{0, q+4}\to X\times\mathbb{R}^{0, q}$ is the projection map,
\item[4)]
$\mathbb{H}\otimes_\mathbb{C}E_i$ is the equivariant Real vector bundles equipped with the Real structure $J\otimes J_{E_i}$,
\item[5)]
$g$ is an equivariant Real vector bundle homomorphism def\/ined by $g(v, w\otimes e)=(v, vw\otimes f(e))$.
\end{enumerate}

One can easily show by generalizing the discussion in the last section of~\cite{AS} that
\begin{Proposition}
\label{shiftbyfour}
$t$ is an isomorphism.
\end{Proposition}

\subsection[The module structure of $KR$-theory of compact simply-connected Lie groups]{The
module structure of $\boldsymbol{KR}$-theory\\ of compact simply-connected Lie groups}
The following structure theorem for $KR$-theory, due to Seymour, is crucial in his computation of
$KR^*(\pt)$-module structure of $KR^*(G)$.
\begin{Theorem}[{\cite[Theorem 4.2]{Se}}]\label{strthm}
Suppose that $K^*(X)$ is a~free abelian group and decomposed by the involution $\overline{\sigma_X^*}$ into the following summands
\begin{gather*}
K^*(X)=M_+\oplus M_-\oplus T\oplus \overline{\sigma_X^*}T,
\end{gather*}
where $\overline{\sigma_X^*}$ is identity on $M_+$ and negation on $M_-$.
Suppose further that there exist $h_1, \dots, h_n\in KR^*(X)$ such that $c(h_1), \dots, c(h_n)$ form a~basis of the
$K^*(+)$-module $K^*(+)\otimes(M_+\oplus M_-)$.
Then, as $KR^*(\pt)$-modules,
\begin{gather*}
KR^*(X)\cong F\oplus r(K^*(+)\otimes T),
\end{gather*}
where $F$ is the free $KR^*(\pt)$-module generated by $h_1, \dots, h_n$.
\end{Theorem}
\begin{Remark}
If $T=0$, then the conditions in Theorem~\ref{strthm} are equivalent to $K^*(X)$ being free abelian and $c: KR^*(X)\to
K^*(X)$ being surjective.
In this special case the theorem
implies that the map $KR^*(X)\otimes_{KR^*(\pt)}K^*(\pt)\to K^*(X)$ def\/ined
by $a\otimes b\mapsto c(a)\cdot b$ is a~ring isomorphism.
This smacks of the def\/inition of weakly equivariant formality (cf.\ Remark~\ref{weakequivformality}) and inspires us to
def\/ine a~similar notion for equivariant $KR$-theory (cf.\ Def\/inition~\ref{realequivformal}).
We say a~real space is \emph{real formal} if it satisf\/ies the conditions of Theorem~\ref{strthm}.
\end{Remark}
\begin{Definition}
Let $\sigma_\mathbb{R}$ be the complex conjugation of $U(n)$ or $U(\infty)$,
and $\sigma_\mathbb{H}$ be the symplectic type involution $g\mapsto J_m\overline{g}J_m^{-1}$ on $U(2m)$, or $U(2\infty)$.
\end{Definition}
For any Real space $X$, $KR^{-1}(X)$ is isomorphic to the abelian group of equivariant homotopy classes of maps from $X$
to $U(\infty)$ which respect $\sigma_X$ and $\sigma_\mathbb{R}$ on $U(\infty)$.
Similarly, $KR^{-5}(X)$, which is isomorphic to $KH^{-1}(X)$ by Proposition~\ref{shiftbyfour}, is isomorphic to the
abelian group of equivariant homotopy classes of maps from $X$ to $U(2\infty)$
which respect $\sigma_X$ and $\sigma_\mathbb{H}$ on $U(2\infty)$
(cf.\ remarks in the last two paragraphs of Appendix of~\cite{Se}).
We can def\/ine maps analogous to those in Def\/inition~\ref{betamap} in the context of $KR$-theory.
\begin{Definition}
Let $\delta_\mathbb{R}: RR(G)\to KR^{-1}(G)$ and $\delta_\mathbb{H}: RH(G)\to KR^{-5}(G)$ be group homomorphisms which
send a~Real (resp.\
Quaternionic) representation to the $KR$-theory element represented by its homotopy class.
\end{Definition}
\begin{Proposition}
If $\rho\in RR(G)$, then $\delta_\mathbb{R}(\rho)$ is represented by the complex of vector bundles in
Proposition~{\rm \ref{-1_class}} equipped with the Real structure given by
\begin{gather*}
\iota: \ G\times\mathbb{R}\times V  \to  G\times\mathbb{R}\times V,
\\
\phantom{\iota: {}} \
(g, t, v)  \mapsto  (\sigma_G(g), t, \overline{v}).
\end{gather*}
If $\rho\in RH(G)$, then $\delta_\mathbb{H}(\rho)$ can be similarly represented, with the Real structure replaced by the
Quaternionic structure.
\end{Proposition}

\emph{From this point on until the end of this section, we further assume that~$G$ is connected and simply-connected unless
otherwise specified}.
It is known that $R(G)$ is a~polynomial ring over $\mathbb{Z}$ generated by fundamental representations, which are
permuted by $\overline{\sigma_{G}^*}$ (cf.\ \cite[Lemma~5.5]{Se}).
Let
\begin{gather*}
R(G)\cong\mathbb{Z}\big[\varphi_1, \dots, \varphi_r, \theta_1, \dots, \theta_s, \gamma_1, \dots, \gamma_t, \overline{\sigma_G^*}
\gamma_1, \dots, \overline{\sigma_G^*}\gamma_t\big],
\end{gather*}
where $\varphi_i\in RR(G, \mathbb{R})$, $\theta_j\in RH(G, \mathbb{R})$, $\gamma_k\in R(G, \mathbb{C})$.
Then $K^*(G)$, as a~free abelian group, is generated by square-free monomials
in $\delta(\varphi_1), {\dots},\delta(\varphi_r)$,
$\delta(\theta_1), {\dots}, \delta(\theta_s)$,
$\delta(\gamma_1), {\dots}, \delta(\gamma_t)$,
$\delta(\overline{\sigma_G^*}\gamma_1)$, ${\dots}, \delta(\overline{\sigma_G^*}\gamma_t)$.
Using Theorem~\ref{strthm}, Seymour obtained

\begin{Theorem}[{\cite[Theorem 5.6]{Se}}]
\label{krtheorymodofg}\quad
\begin{enumerate}\itemsep=0pt
\item[$1.$]
Suppose that $\overline{\sigma_G^*}$ acts as identity on $R(G)$, i.e.\
any irreducible Real representation of~$G$ is either of real type or quaternionic type.
Then as $KR^*(\pt)$-modules,
\begin{gather*}
KR^*(G)\cong\wedge_{KR^*(\pt)}(\delta_\mathbb{R}(\varphi_1), \dots, \delta_\mathbb{R}(\varphi_r),
\delta_\mathbb{H}(\theta_1), \dots, \delta_\mathbb{H}(\theta_s)).
\end{gather*}
\item[$2.$]
More generally, $c(\delta_\mathbb{R}(\varphi_i))=\delta(\varphi_i)$,
$c(\delta_\mathbb{H}(\theta_j))=\beta^2\cdot\delta(\theta_j)$, and there exist $\lambda_1, \dots, \lambda_t\in KR^0(G)$
such that $c(\lambda_k)=\beta^3\cdot\delta(\gamma_k)\delta(\overline{\sigma_G^*}
\gamma_k)$, and
\begin{gather*}
KR^*(G)\cong P\oplus T\cdot P
\end{gather*}
as $KR^*(\pt)$-module, where
\begin{itemize}\itemsep=0pt
\item
$P\cong\bigwedge_{KR^*(\pt)}(\delta_\mathbb{R}(\varphi_1), \dots, \delta_\mathbb{R}(\varphi_r),
\delta_\mathbb{H}(\theta_1), \dots, \delta_\mathbb{H}(\theta_s), \lambda_1, \dots, \lambda_t)$,
\item
$T$ is the additive abelian group generated by the set
\begin{gather*}
\{r(\beta^i\cdot\delta(\gamma_1)^{\varepsilon_1}\cdots\delta(\gamma_t)^{\varepsilon_t}\delta(\overline{\sigma_G^*}
\gamma_1)^{\nu_1}\cdots\delta(\overline{\sigma_G^*}\gamma_t)^{\nu_t})\},
\end{gather*}
where $\varepsilon_1, \dots, \varepsilon_t$, $\nu_1, \dots, \nu_t$ are either $0$ or $1$, $\varepsilon_k$ and~$\nu_k$ are not
equal to 1 at the same time for $1\leq k\leq t$, and the first index $k_0$ where $\varepsilon_{k_0}=1$ is less than the
first index~$k_1$ where $\nu_{k_1}=1$.
\end{itemize}
Moreover,
\begin{enumerate}\itemsep=0pt
\item[$(a)$]
$\lambda_k^2=0$ for all $1\leq k\leq t$,
\item[$(b)$]
$\delta_\mathbb{R}(\varphi_i)^2$ and $\delta_\mathbb{H}(\theta_j)^2$ are divisible by $\eta$.
\end{enumerate}
\end{enumerate}
\end{Theorem}

\begin{Definition}
Let $\omega_t:=\delta_{\epsilon_t, 1-\nu_t}$ and
\begin{gather*}
r_{i, \varepsilon_1, \dots, \varepsilon_t, \nu_1, \dots,
\nu_t}:=r\left(\beta^i\cdot\delta(\gamma_1)^{\varepsilon_1}\cdots\delta(\gamma_t)^{\varepsilon_t}\delta(\overline{\sigma_G^*}
(\gamma_1))^{\nu_1}\cdots\delta\big(\overline{\sigma_G^*}(\gamma_t)\big)^{\nu_t}\right)\in T.
\end{gather*}
\end{Definition}
\begin{Corollary}\label{krglooseend}\quad
\begin{enumerate}\itemsep=0pt
\item[$1.$]
$KR^*(G)$ is generated by $\delta_\mathbb{R}(\varphi_1),  {\dots}, \delta_\mathbb{R}(\varphi_r)$,
$\delta_\mathbb{H}(\theta_1),  {\dots}, \delta_\mathbb{H}(\theta_s)$, $\lambda_1,  {\dots}, \lambda_t\!$ and $r_{i,
\varepsilon_1, {\dots}, \varepsilon_t, \nu_1, {\dots}, \nu_t}$ $\in T$ as an algebra over $KR^*(\pt)$.
\item[$2.$]
\begin{gather*}
r_{i, \varepsilon_1, \dots, \varepsilon_t, \nu_1, \dots, \nu_t}^2=
\begin{cases}
\eta^2\lambda_1^{\omega_1}\cdots\lambda_t^{\omega_t},
&\text{if }r_{i, \varepsilon_1, \dots, \varepsilon_t, \nu_1, \dots, \nu_t}\text{ is of degree }-1\text{ or }-5,
\\
\pm\mu\lambda_1^{\omega_1}\cdots\lambda_t^{\omega_t},
&\text{if }r_{i, \varepsilon_1, \dots, \varepsilon_t, \nu_1, \dots, \nu_t}\text{ is of degree }-2\text{ or }-6,
\\
0
&\text{otherwise}.
\end{cases}
\end{gather*}
The sign depends on $i$, $\varepsilon_1, \dots, \varepsilon_t$, $\nu_1, \dots, \nu_t$ and can be determined using formulae
from~$(2)$ of Proposition~{\rm \ref{krprelim}}.
\item[$3.$]
$r_{i, \varepsilon_1, \dots, \varepsilon_t, \nu_1, \dots, \nu_t}\eta=0$, and $r_{i, \varepsilon_1, \dots, \varepsilon_t,
\nu_1, \dots, \nu_t}\mu=2r_{i+2, \varepsilon_1, \dots, \varepsilon_t, \nu_1, \dots, \nu_t}$.
\end{enumerate}
\end{Corollary}

\begin{proof}
The Corollary follows easily from the various properties of the realif\/ication map and the complexif\/ication map in
Proposition~\ref{krprelim}, and the fact that $c(\lambda_k)=\beta^3\cdot\delta(\gamma_k)\delta(\overline{\sigma_G^*}
\gamma_k)$.
For example,
\begin{gather*}
r_{i, \varepsilon_1, \dots, \varepsilon_t, \nu_1, \dots,
\nu_t}\eta=r\big(\beta^i\cdot\delta(\gamma_1)^{\varepsilon_1}\cdots\delta(\gamma_t)^{\varepsilon_t}\delta(\overline{\sigma_G^*}
(\gamma_1))^{\nu_1}\cdots\delta\big(\overline{\sigma_G^*}(\gamma_t)\big)^{\nu_t}c(\eta)\big)
=0,
\\
r_{i, \varepsilon_1, \dots, \varepsilon_t, \nu_1, \dots,\nu_t}\mu
=r\big(\beta^i\cdot\delta(\gamma_1)^{\varepsilon_1}\cdots\delta(\gamma_t)^{\varepsilon_t}\delta(\overline{\sigma_G^*}
(\gamma_1))^{\nu_1}\cdots\delta\big(\overline{\sigma_G^*}(\gamma_t)\big)^{\nu_t}c(\mu)\big)
\\
\phantom{r_{i, \varepsilon_1, \dots, \varepsilon_t, \nu_1, \dots,\nu_t}\mu}
=r_{i+2, \varepsilon_1, \dots, \varepsilon_t, \nu_1, \dots, \nu_t}.\tag*{\qed}
\end{gather*}
\renewcommand{\qed}{}
\end{proof}

In fact Theorems~\ref{hodgkin} and~\ref{strthm} also yield the following description of module structure of
$KR$-theory of a~compact connected Real Lie group with torsion-free fundamental group with a~restriction on the types of
the Real representations.
\begin{Theorem}
\label{spkrtheorymodofg}
Let~$G$ be a~compact connected real Lie group with $\pi_1(G)$ torsion-free.
Suppose that $R(G, \mathbb{C})=0$, i.e.\
$\overline{\sigma_G^*}$ acts as identity on $R(G)$.
Then $KR^*(G)$ is isomorphic to $\wedge_{KR^*(\pt)}^*(\text{\rm Im}(\widetilde{\delta}_\mathbb{R}),
\text{\rm Im}(\widetilde{\delta}_\mathbb{H}))$ as $KR^*(\pt)$-modules.
\end{Theorem}

As we see from Theorem~\ref{krtheorymodofg} and Corollary~\ref{krglooseend}, to get a~full description of the ring
structure of $KR^*(G)$, it remains to f\/igure out $\delta_\mathbb{R}(\varphi_i)^2$ and $\delta_\mathbb{H}(\theta_j)^2$.
We will, in the end, obtain formulae for the squares by way of computing the ring structure of $KR_G^*(G)$ and applying
the forgetful map.
In particular, we will show that $\delta_\mathbb{R}(\varphi_i)^2$ and $\delta_\mathbb{H}(\theta_j)^2$ in general are
non-zero.
So, unlike the complex $K$-theory, $KR^*(G)$ is not an exterior algebra in general.
Nevertheless, $KR^*(G)$ is not far from being an exterior algebra, in the sense of the following
\begin{Corollary}\quad
\begin{enumerate}\itemsep=0pt
\item[$1.$]
$KR^*(\pt)_2$, which is the ring obtained by inverting the prime~$2$ in $KR^*(\pt)$, is isomorphic to
$\mathbb{Z}\left[\frac{1}{2}, \mu\right]/(\mu^2-4)\cong\mathbb{Z}\left[\frac{1}{2}, \beta^2\right]/((\beta^2)^2-1)$.
\item[$2.$]
Suppose that $R(G, \mathbb{C})=0$.
$KR^*(G)_2$, which is the ring obtained by inverting the prime~$2$ in $KR^*(G)$, is isomorphic to, as
$KR^*(\pt)_2$-algebra
\begin{gather*}
\bigwedge\nolimits_{KR^*(\pt)_2}\big(\delta_\mathbb{R}(\varphi_1), \dots, \delta_\mathbb{R}(\varphi_r),
\delta_\mathbb{H}(\theta_1), \dots, \delta_\mathbb{H}(\theta_s)\big).
\end{gather*}
\end{enumerate}
\end{Corollary}

\section[The coef\/f\/icient ring $KR_G^*(\pt)$]{The coef\/f\/icient ring $\boldsymbol{KR_G^*(\pt)}$}\label{Section3}

In this section, we assume that~$G$ is a~compact Real Lie group, and will prove a~result on the coef\/f\/icient ring
$KR_G^*(\pt)$.
In~\cite{AS}, all graded pieces of $KR_G^*(\pt)$ were worked out using Real Clif\/ford~$G$-modules.
We record them in the following
\begin{Proposition}
\label{ASKRcomp}
$KR^{-q}_G(\pt)$, as abelian groups, for $0\leq q\leq 7$, are isomorphic to $RR(G)$, $RR(G)/\rho(R(G))$,
$R(G)/j(RH(G))$, $0$, $RH(G)$, $RH(G)/\eta(R(G))$, $R(G)/i(RR(G))$ and $0$ respectively, where the maps $i$, $j$, $\rho$, $\eta$
are as in Propositions~{\rm \ref{realforgetfulinj}} and~{\rm \ref{quatforgetfulinj}}.
\end{Proposition}
\begin{Remark}
Note from the above proposition
that $KR_G^0(\pt)\oplus KR^{-4}_G(\pt)\cong RR(G)\oplus RH(G)$.
In this way we can view $RR(G)\oplus RH(G)$ as a~graded ring where $RR(G)$ is of degree 0 and $RH(G)$ of degree $-4$.
\end{Remark}
\begin{Proposition}\label{KR_G}\quad
\begin{enumerate}\itemsep=0pt
\item[$1.$]
Suppose $R(G, \mathbb{C})=0$.
Then the map
\begin{gather*}
f:  \ (RR(G, \mathbb{R})\oplus RH(G, \mathbb{R}))\otimes KR^*(\pt)  \to  KR^*_G(\pt),
\\
\phantom{f: {}} \
\rho_1\otimes x_1\oplus \rho_2\otimes x_2  \mapsto  \rho_1\cdot x_1+\rho_2\cdot x_2
\end{gather*}
is an isomorphism of graded rings.
\item[$2.$]
In general,
\begin{gather*}
f: \ (RR(G, \mathbb{R})\oplus RH(G, \mathbb{R}))\otimes KR^*(\pt)\oplus r(R(G, \mathbb{C})\otimes K^*(+))\to
KR^*_G(\pt),
\\
\phantom{f: {}} \
\rho_1\otimes x_1\oplus \rho_2\otimes x_2\oplus r(\rho_3\otimes \beta^i)\mapsto \rho_1\cdot x_1+\rho_2\cdot
x_2+r(\rho_3\cdot\beta^i)
\end{gather*}
is an isomorphism of graded abelian groups.
\item[$3.$]
If $\rho$ is an irreducible complex representation of complex type, then $\eta r(\beta^i\cdot\rho)=0$ and $\mu
r(\beta^i\cdot\rho)=2r(\beta^{i+2}\cdot\rho)$.
\end{enumerate}
\end{Proposition}

\begin{proof}
The proposition
follows by verifying the isomorphism in dif\/ferent degree pieces against the description in
Proposition~\ref{ASKRcomp}.
For example, in degree 0,
\begin{gather*}
RR(G, \mathbb{R})\otimes KR^0(\pt)\oplus RH(G, \mathbb{R})\otimes KR^{-4}(\pt)\oplus r(R(G,
\mathbb{C})\otimes K^0(+))
\\
\qquad=RR(G, \mathbb{R})\oplus RH(G, \mathbb{R})\otimes \mathbb{Z}\mu\oplus RR(G, \mathbb{C})
\\
\qquad
\cong RR(G, \mathbb{R})\oplus RR(G, \mathbb{H})\oplus RR(G, \mathbb{C})
\\
(\text{if} \ \ [V]\in RH(G, \mathbb{R}),\ \ \text{then} \ \ [V]\cdot\mu=[V\oplus V]\in RR(G, \mathbb{H}))
\\
\qquad=RR(G)
=KR^0_G(\pt).
\end{gather*}
(3) follows from Proposition~\ref{krprelim}.
\end{proof}
\begin{Remark}
In~\cite{AS}, $KR_G^*(X)$, where the~$G$-action is trivial, is given as the following direct sum of abelian groups
\begin{gather*}
RR(G, \mathbb{R})\otimes KR^*(X)\oplus RR(G, \mathbb{C})\otimes KC^*(X)\oplus RR(G, \mathbb{H})\otimes KH^*(X),
\end{gather*}
where $KC^*(X)$ and $KH^*(X)$ are Grothendieck groups of the so-called `Complex vector bundles' and `Quaternionic vector
bundles' of $X$.
We f\/ind Proposition~\ref{KR_G}, which is motivated by this description, better because the ring structure of the
coef\/f\/icient ring is more apparent when cast in this light.
The proposition
is, as we will see in the next section, a~consequence of a~structure theorem of equivariant $KR$-theory
(Theorem~\ref{equivstrthm}), and therefore still holds true if the point is replaced by any general space $X$ with
trivial~$G$-action.
\end{Remark}

\section[Equivariant $KR$-theory rings of compact simply-connected Lie groups]{Equivariant
$\boldsymbol{KR}$-theory rings \\ of compact simply-connected Lie groups}\label{Section4}

Throughout this section we assume that~$G$ is a~compact, connected and simply-connected Real Lie group unless otherwise
specif\/ied.
We will prove the main result of this paper, Theorem~\ref{mainthm2}, which gives the ring structure of $KR_G^*(G)$.
Our strategy is outlined as follows.
\begin{enumerate}\itemsep=0pt
\item
We obtain a~result on the structure of $KR^*_G(G)$ (Corollary~\ref{equivstrthmg}) which is analogous to
Theorem~\ref{krtheorymodofg} and Proposition~\ref{KR_G}.
We def\/ine $\delta_\mathbb{R}^G(\varphi_i)$, $\delta_\mathbb{H}^G(\theta_j)$, $\lambda_k^G$ and $r^G_{\rho, i,
\varepsilon_1, \dots, \varepsilon_t, \nu_1, \dots, \nu_t}$ (cf.\ Def\/inition~\ref{equivlift} and
Corollary~\ref{eqmodgenerators}), which generate $KR_G^*(G)$ as a~$KR_G^*(\pt)$-algebra, as a~result of
Corollary~\ref{equivstrthmg}.
We show that $(\lambda_k^G)^2=0$ (cf.\ Proposition~\ref{cplxsquare}).
\item
We compute the module structure of $KR^*_{(U(n),
\sigma_\mathbb{F})}(U(n),
\sigma_\mathbb{F})$ for $\mathbb{F}=\mathbb{R}$ and $\mathbb{H}$.
\item
Let $T$ be the maximal torus of diagonal matrices in $U(n)$ and, by abuse of notation, $\sigma_\mathbb{R}$ be the inversion map on $T$, $\sigma_\mathbb{H}$ be the involution on $U(n)/T$ (where $n=2m$ is even) def\/ined by $gT\mapsto J_m\overline{g}T$.
We show that the restriction map
\begin{gather*}
p^*_G: \  KR^*_{(U(n),
\sigma_\mathbb{R})}(U(n),
\sigma_\mathbb{R})\to KR^*_{(T,
\sigma_\mathbb{R})}(T,
\sigma_\mathbb{R})
\end{gather*}
and the map
\begin{gather*}
q^*_G: \ KR^*_{(U(2m),
\sigma_\mathbb{H})}(U(2m),
\sigma_\mathbb{H})\to KR^*_{(U(2m),
\sigma_\mathbb{H})}(U(2m)/T\times T,
\sigma_\mathbb{H}
\times
\sigma_\mathbb{R})
\end{gather*}
induced by the Weyl covering map $q_G: U(2m)/T\times T\to U(2m)$, $(gT, t)\mapsto gtg^{-1}$, are injective.
\item
Let $\sigma_n$ be the class of the standard representation of $U(n)$.
We pass the computation of the two squares $\delta_\mathbb{R}^G(\sigma_n)^2\in KR^*_{(U(n), \sigma_\mathbb{R})}(U(n),
\sigma_\mathbb{R})$ and $\delta_\mathbb{H}^G(\sigma_{2m})^2\in KR^*_{(U(2m),
\sigma_\mathbb{H})}(U(2m),
\sigma_\mathbb{H})$ through the induced map $p^*_G$ and $q_G^*$ to their images and get equations~\eqref{realsigmasquare}
and~\eqref{quatsigmasquare} in Proposition~\ref{equivkrthysquaresigma}.
\item
Applying induced maps $\varphi_i^*$ and $\theta_j^*$ to equations~\eqref{realsigmasquare} and~\eqref{quatsigmasquare}
yields equations~\eqref{realsquare} and~\eqref{quatsquare} in Theorem~\ref{equivkrthysquare} which, together with
Proposition~\ref{cplxsquare} and some relations among~$\eta$,~$\mu$ and $r^G_{\rho, i, \varepsilon_1, \dots,
\varepsilon_t, \nu_1, \dots, \nu_t}$ deduced from Proposition~\ref{krprelim}, describe completely the ring structure of
$KR_G^*(G)$ (cf.\ Theorem~\ref{mainthm2}).
\end{enumerate}
\begin{Remark}\quad
\begin{enumerate}\itemsep=0pt
\item
Seymour f\/irst suggested the analogues of Steps 3, 4 and 5 in the ordinary $KR$-theory case in~\cite{Se} in an attempt to
compute $\delta_\mathbb{R}(\varphi_i)^2$ and $\delta_\mathbb{H}(\theta_j)^2$, but failed to establish Step~3, which he
assumed to be true to make conjectures about $\delta_\mathbb{R}(\varphi_i)^2$.
\item
In equivariant complex $K$-theory, $K_G^*(G/T\times T)\cong K_T^*(T)$ for any compact Lie group~$G$, and the two maps
$p_G^*$ (the restriction map induced by the inclusion $T\hookrightarrow G$) and $q_G^*$ which is induced by the Weyl
covering map are the same.
If $\pi_1(G)$ is torsion-free, then these two maps are shown to be injective (cf.~\cite{BZ}.
In fact it is even shown there that the maps inject onto the Weyl invariants of $K_T^*(T)$).
In the case of equivariant $KR$-theory, things are more complicated.
First of all, while in the case where $(G,\sigma_G)=(U(n), \sigma_\mathbb{R})$, it is true that $KR_G^*(G/T\times T)\cong KR_T^*(T)$, and $p_G^*$ and $q_G^*$ are the same, it is no longer true in
the case where $(G,\sigma_G)=(U(2m), \sigma_\mathbb{H})$.
In Step 3, we use $q_G^*$ for the quaternionic type involution case because we f\/ind that it admits an easier description
than $p_G^*$ does.
Second, we do not know whether $p_G^*$ and $q_G^*$ are injective for general compact Real Lie groups (equipped with any
Lie group involution).
For our purpose it is suf\/f\/icient to show the injectivity results in Step 3.
\end{enumerate}
\end{Remark}

\subsection{A structure theorem}
\begin{Definition}
\label{realequivformal}
A~$G$-space $X$ is a~\emph{Real equivariantly formal} space if
\begin{enumerate}\itemsep=0pt
\item[1)]
$G$ is a~compact Real Lie group,
\item[2)]
$X$ is a~weakly equivariantly formal~$G$-space, and
\item[3)]
the forgetful map $KR^*_G(X)\to KR^*(X)$ admits a~section $s_R: KR^*(X)\to KR^*_G(X)$ which is
a~$KR^*(\pt)$-module homomorphism.
\end{enumerate}
\end{Definition}
\begin{Remark}
If $X$ is a~weakly equivariantly formal~$G$-space, then the forgetful map $K_G^*(X)\to K^*(X)$ admits a~(not necessarily
unique) section $s: K^*(X)\to K_G^*(X)$ which is a~group homomorphism.
\end{Remark}
\begin{Definition}
For a~section $s: K^*(X)\to K^*_G(X)$ (resp.\
$s_R: KR^*(X)\to KR_G^*(X))$ and $a\in K^*(X)$ (resp.\
$a\in KR^*(X)$), we call $s(a)$ (resp.~$s_R(a)$) a~\emph{$($Real$)$ equivariant lift} of~$a$, with respect to $s$ (resp.~$s_R$).
\end{Definition}

We f\/irst prove a~structure theorem of equivariant $KR$-theory of Real equivariantly formal spaces.
\begin{Theorem}
\label{equivstrthm}
Let $X$ be a~Real equivariantly formal space.
For any element $a\in K^*(X)$ $($resp.\
$a\in KR^*(X))$, let $a_G\in K_G^*(X)$ $($resp.\
$a_G\in KR_G^*(X))$ be a~$($Real$)$ equivariant lift of $a$ with respect to a~group homomorphic section $s$ $($resp.~$s_R$ which is a~$KR^*(\pt)$-module homomorphism$)$.
Then the map
\begin{gather*}
\begin{split}
& f:  \ (RR(G, \mathbb{R})\oplus RH(G, \mathbb{R}))\otimes KR^*(X)\oplus r(R(G, \mathbb{C})\otimes K^*(X))  \to  KR_G^*(X),
\\
& \phantom{f: {}} \
\rho_1\otimes a_1\oplus r(\rho_2\otimes a_2)  \mapsto  \rho_1\cdot (a_1)_{G}\oplus r(\rho_2\cdot (a_2)_G).
\end{split}
\end{gather*}
is a~group isomorphism.
In particular, if $R(G, \mathbb{C})=0$, then $f$ is a~$KR_G^*(\pt)$-module isomorphism.
\end{Theorem}

\begin{proof}
Consider the following $H(p, q)$-systems
\begin{gather*}
HR^\alpha(p, q):=KR^{-\alpha}\big(X\times S^{q, 0}, X\times S^{p, 0}\big)\cong KR^{-\alpha+p}\big(X\times S^{q-p, 0}\big),
\\
H^{\alpha}(p, q):= K^{-\alpha}\big(X\times S^{q, 0}, X\times S^{p, 0}\big)\cong K^{-\alpha+p}\big(X\times S^{q-p, 0}\big),
\\
HR_G^{\alpha}(p, q):=KR_G^{-\alpha}\big(X\times S^{q, 0}, X\times S^{p, 0}\big)\cong KR_G^{-\alpha+p}\big(X\times S^{q-p, 0}\big).
\end{gather*}
For the last $H(p, q)$-system,~$G$ acts on $S^{q-p, 0}$ trivially.
The spectral sequences induced by these $H(p, q)$-systems converge to $KR^*(X)$, $K^*(X)$ and $KR_G^*(X)$ respectively
(for the assertion for the f\/irst two $H(p, q)$-systems,
see the proofs of Theorem~3.1 and Lemma~4.1 of~\cite{Se}.
That the third $H(p, q)$-system converges to $KR_G^*(X)$ follows from a~straightforward generalization of the
aforementioned proofs by adding equivariant structure throughout).
Consider the two long exact sequences for the pair $(X\times B^{q-p, 0}, X\times S^{q-p, 0})$, with the top exact
sequence involving equivariant $KR$-theory and the bottom one ordinary $KR$-theory, and the vertical maps being
forgetful maps.
By applying the f\/ive-lemma, we have that each element in the f\/irst two $H(p, q)$-systems has a~(Real) equivariant lift.
Def\/ine a~group homomorphism
\begin{gather*}
f(p, q): \ (RR(G, \mathbb{R})\oplus RH(G, \mathbb{R}))\otimes HR^\alpha(p, q)\oplus r(R(G, \mathbb{C})\otimes H^\alpha(p,
q))\to HR_G^{\alpha}(p, q)
\end{gather*}
by
\begin{gather*}
\rho_1\otimes a_1\oplus r(\rho_2\otimes a_2)\mapsto \rho_1\cdot (a_1)_{G}\oplus r(\rho_2\cdot(a_2)_G).
\end{gather*}
As $RR(G, \mathbb{R})$, $RH(G, \mathbb{R})$ and $R(G, \mathbb{C})$ are free abelian groups, and tensoring free abelian
groups and taking cohomology commute, $f$ is the abutment of $f(p, q)$.
On the $E_1^{p, q}$-page, $f(p, q)$ becomes
\begin{gather*}
f_1^{p, q}: \ (RR(G, \mathbb{R})\oplus RH(G, \mathbb{R}))\otimes KR^{-q}\big(X\times S^{1, 0}\big)\oplus r(R(G, \mathbb{C})\otimes
K^{-q}\big(X\times S^{1, 0}\big))
\\
\hphantom{f_1^{p, q}: {}}{} \ \to
KR_G^{-q}\big(X\times S^{1, 0}\big).
\end{gather*}
Note that $K^{-q}(X\times S^{1, 0})\cong K^{-q}(X)\oplus K^{-q}(X)$, $KR^{-q}(X\times S^{1, 0})\cong K^{-q}(X)$, and
$KR_G^{-q}(X\times S^{1, 0})\cong K_G^{-q}(X)$.
With the above identif\/ication,
\begin{gather*}
r: \ R(G, \mathbb{C})\otimes K^{-q}\big(X\times S^{1, 0}\big)  \to  KR_G^{-q}\big(X\times S^{1, 0}\big)\cong K_G^{-q}(X),
\\
\phantom{r: {}} \
\rho_1\otimes(a_1, 0)\oplus \rho_2\otimes (0, a_2)  \mapsto  \rho_1\cdot (a_1)_G+\overline{\sigma_G^*}
\rho_2\cdot \big(\overline{\sigma_G^*}a_2\big)_G.
\end{gather*}
So $r(R(G, \mathbb{C})\otimes K^{-q}(X\times S^{1, 0}))= R(G, \mathbb{C})\otimes K^{-q}(X)$.
$f_1^{p, q}$ is a~group homomorphism from $(RR(G, \mathbb{R})\oplus RH(G, \mathbb{R})\oplus R(G, \mathbb{C}))\otimes
K^{-q}(X)\cong R(G)\otimes K^{-q}(X)$ to $K_G^{-q}(X)$, which is an isomorphism by weak equivariant formality of $X$.
It follows that $f$ is also an isomorphism.
If $R(G, \mathbb{C})=0$, then by~(1) of Proposition~\ref{KR_G} and~(3) of Def\/inition~\ref{realequivformal}, $f$ is
indeed a~$KR_G^*(\pt)$-module isomorphism.
\end{proof}
\begin{Remark}
The term `Real equivariant formality' is suggested by the observation that, if $X$ is a~Real equivariantly formal
$G$-space and $R(G, \mathbb{C})=0$, then the map
\begin{gather*}
KR^*_G(X)\otimes_{RR(G, \mathbb{R})\oplus RH(G, \mathbb{R})}\mathbb{Z}\to KR^*(X)
\end{gather*}
induced by the forgetful map is a~ring isomorphism, which smacks of the ring isomorphism in the def\/inition of weak
equivariant formality.
\end{Remark}
\begin{Lemma}
\label{realequivlift}
$\delta_\mathbb{R}(\varphi_i)$, $\delta_\mathbb{H}(\theta_j)$, $\lambda_k$ and $r_{i, \varepsilon_1, \dots,
\varepsilon_t, \nu_1, \dots, \nu_t}\in KR^*(G)$ all have Real equivariant lifts in $KR_G^*(G)$.
Hence~$G$ is a~Real equivariantly formal space.
\end{Lemma}
\begin{proof}
A natural choice of a~Real equivariant lift of $\delta_\mathbb{R}(\varphi_i)$ is represented by the complex of vector
bundles in Proposition~\ref{-1_class} equipped with both the Real structure and equivariant structure def\/ined for
$\delta_\mathbb{R}(\varphi_i)$ and $\delta_G(\varphi)$ respectively.
These two structures are easily seen to be compatible.
A Real equivariant lift of $\delta_\mathbb{H}(\theta_j)$ can be similarly def\/ined.
The class
\begin{gather*}
r\big(\beta^i\cdot\delta_G(\gamma_1)^{\varepsilon_1}\cdots\delta_G(\gamma_t)^{\varepsilon_t}\delta_G(\overline{\sigma_G^*}
\gamma_1)^{\nu_1}\cdots\delta_G(\overline{\sigma_G^*}
\gamma_t)^{\nu_t}\big)
\end{gather*}
obviously is a~Real equivariant lift of $r_{i, \varepsilon_1, \dots, \varepsilon_t, \nu_1, \dots, \nu_t}$.
By adding the natural equivariant structure throughout the construction of $\lambda_k$ in the proof of Proposition~4.6
in~\cite{Se}, one can obtain a~Real equivariant lift of $\lambda_k$.
\end{proof}

\begin{Definition}
\label{equivlift}
We f\/ix a~choice of equivariant lift of any element $a\in K^*(G)$ by def\/ining $\delta_G(\rho)$ to be the equivariant lift
of $\delta(\rho)$.
Similarly, we f\/ix a~choice of Real equivariant lift of $a\in KR^*(G)$ by def\/ining $\delta_\mathbb{R}^G(\varphi_i)$,
$\delta_\mathbb{H}^G(\theta_j)$, $\lambda_k^G$, and $r^G_{i, \varepsilon_1, \dots, \varepsilon_t, \nu_1, \dots, \nu_t}$
in the proof of Lemma~\ref{realequivlift} to be the Real equivariant lift of $\delta_\mathbb{R}(\varphi_i)$,
$\delta_\mathbb{H}(\theta_j)$, $\lambda_k$ and $r_{i, \varepsilon_1, \dots, \varepsilon_t, \nu_1, \dots, \nu_t}$.
\end{Definition}
\begin{Remark}
\label{cplxeqsquare}
$\lambda^G_k$ satisf\/ies $c(\lambda^G_k)=\beta^3\delta_G(\gamma_k)\delta_G(\overline{\sigma_G^*}
\gamma_k)$.
\end{Remark}
\begin{Corollary}
\label{equivstrthmg}
Let~$G$ be a~compact, connected and simply-connected Real Lie group.
The map
\begin{gather*}
f: \ (RR(G, \mathbb{R})\oplus RH(G, \mathbb{R}))\otimes KR^*(G)\oplus r(R(G, \mathbb{C})\otimes K^*(G))  \to  KR_G^*(G),
\\
\phantom{f: {}} \
\rho_1\otimes a_1\oplus r(\rho_2\otimes a_2)  \mapsto  \rho_1\cdot (a_1)_{G}\oplus r(\rho_2\cdot (a_2)_G)
\end{gather*}
is a~group isomorphism.
Here $(a_i)_G$ is the $($Real$)$ equivariant lift defined as in Definition~{\rm \ref{equivlift}}.
In particular, if $R(G, \mathbb{C})=0$, then $f$ is an isomorphism of $KR^*_{G}(\pt)$-modules from
$KR_G^*(\pt)\otimes K^*(G)$ to $KR_G^*(G)$.
\end{Corollary}
\begin{proof}
The result follows from Theorem~\ref{equivstrthm} and Lemma~\ref{realequivlift}.
In the special case where $R(G, \mathbb{C})=0$, $KR^*(G)$ is isomorphic to $KR^*(\pt)\otimes K^*(G)$ as
$KR^*(\pt)$-modules by (1) of Theorem~\ref{krtheorymodofg}, and applying Theorem~\ref{equivstrthm} and
Proposition~\ref{KR_G} give $KR^*_G(G)\cong RR(G)\otimes KR^*(G)\cong RR(G)\otimes KR^*(\pt)\otimes K^*(G)\cong
KR^*_G(\pt)\otimes K^*(G)$.
In this way $f$ is a~$KR_G^*(\pt)$-module isomorphism from $KR_G^*(\pt)\otimes K^*(G)$ to $KR_G^*(G)$.
\end{proof}
\begin{Corollary}
\label{eqmodgenerators}
Let
\begin{gather*}
r^G_{\rho, i, \varepsilon_1, \dots, \varepsilon_t, \nu_1, \dots,\nu_t}
:=r\big(\beta^i\cdot\rho\delta_G(\gamma_1)^{\varepsilon_1}\cdots\delta_G(\gamma_t)^{\varepsilon_t}\delta_G(\overline{\sigma_G^*}
\gamma_1)^{\nu_1}\cdots\delta_G(\overline{\sigma_G^*}
\gamma_t)^{\nu_t}\big),
\end{gather*}
where $\rho\in R(G, \mathbb{C})\oplus \mathbb{Z}\cdot\rho_{\text{\rm triv}}$ and $\varepsilon_1, \dots, \varepsilon_t$,
$\nu_1, \dots, \nu_t$ are as in Theorem~{\rm \ref{strthm}}.
Then $KR_G^*(G)$, as an algebra over $KR_G^*(\pt)$, is generated by $\delta_\mathbb{R}^G(\varphi_1), \dots,
\delta_\mathbb{R}^G(\varphi_r)$, $\delta_\mathbb{H}^G(\theta_1), \dots, \delta_\mathbb{H}^G(\theta_s)$, $\lambda_1, \dots,
\lambda_t$, and $r^G_{\rho, i, \varepsilon_1, \dots, \varepsilon_t, \nu_1, \dots, \nu_t}$.
\end{Corollary}
\begin{Remark}
If $\rho=\rho_{\text{triv}}$, then $r^G_{\rho, i, \varepsilon_1, \dots, \varepsilon_t, \nu_1, \dots, \nu_t}=r^G_{i,
\varepsilon_1, \dots, \varepsilon_t, \nu_1, \dots, \nu_t}$.
If $\rho\in R(G, \mathbb{C})$, then $r^G_{\rho, i, \varepsilon_1, \dots, \varepsilon_t, \nu_1, \dots, \nu_t}$ comes from
$r(R(G, \mathbb{C})\otimes K^*(G))$ in the decomposition of Theorem~\ref{equivstrthm}.
\end{Remark}

Now we are in a~position to compute $(\lambda_k^G)^2$ by imitating the proof of Proposition~4.7 in~\cite{Se}.
\begin{Proposition}
\label{cplxsquare}
$\big(\lambda_k^G\big)^2=0$.
\end{Proposition}
\begin{proof}
Consider the Real Lie group $(U(n)\times U(n),
\sigma_\mathbb{C})$, where $\sigma_\mathbb{C}
(g_1, g_2)=(\overline{g_2}, \overline{g_1})$.
Let $p_j: U(n)\times U(n)\to U(n)$ be the projection onto the $j$-th factor, and $u_i=p_1^*(\wedge^i
\sigma_n)$, $v_i=p_2^*(\wedge^i\sigma_n)$.
Thus $\overline{\sigma_\mathbb{C}^*}u_i=v_i$.
A decomposition of $K^*(U(n)\times U(n))$ by the induced involution $\overline{\sigma_\mathbb{C}^*}$ is given by $M\oplus T\oplus \overline{\sigma_\mathbb{C}^*}T$, where $M$ is the subalgebra generated by $\delta(u_1)\delta(v_1), \dots, \delta(u_n)\delta(v_n)$.
By Proposition~\ref{strthm}, there exist $h_1, \dots, h_n\in KR^0(U(n)\times U(n),
\sigma_\mathbb{C})$ such that $c(h_i)=\beta^3\delta(u_i)\delta(v_i)$, and $KR^*(U(n)\times U(n),
\sigma_\mathbb{C})\cong F\oplus r(K^*(+)\otimes T)$, where $F$ is the $KR^*(\pt)$-module freely generated by monomials in $h_1,
\dots, h_n$.
By Corollary~\ref{equivstrthmg},
\begin{gather*}
KR^*_{(U(n)\times U(n),
\sigma_\mathbb{C})}(U(n)\times U(n),
\sigma_\mathbb{C})
\cong
RR(U(n)\times U(n),
\sigma_\mathbb{C}, \mathbb{R})\otimes(F\oplus r(K^*(+)\otimes T))
\\
\qquad{}
\oplus
r(R(U(n)\times U(n),
\sigma_\mathbb{C}, \mathbb{C})\otimes K^*(U(n)\times U(n))).
\end{gather*}
Let $h_i^G$ be the equivariant lift of $h_i$ as def\/ined in Def\/inition~\ref{equivlift}.
So $c(h_i^G)=\beta^3\delta_G(u_i)\delta_G(v_i)$ and $c((h_i^G)^2)=0$.
Consequently $(h_i^G)^2=\eta k_i$ for some $k_i\in KR^{-7}_{(U(n)\times U(n),
\sigma_\mathbb{C})}(U(n)\times U(n),
\sigma_\mathbb{C})$ (cf.\ Gysin sequence~(3.4) in~\cite{At3} and its equivariant analogue).
Since $\eta\cdot r(\cdot)=0$, we may assume that~$k_i$ is from the component $RR(U(n)\times U(n),
\sigma_\mathbb{C}, \mathbb{R})\otimes F$.
But the degree $-7$ piece of the later is~0.
So $(h_i^G)^2=0$.

Consider the map
\begin{gather*}
\gamma_k\times \overline{\sigma_\mathbb{C}^*}\gamma_k: \ (G,\sigma_G)  \to (U(n)\times U(n), \sigma_\mathbb{C}),
\\
\phantom{\gamma_k\times \overline{\sigma_\mathbb{C}^*}\gamma_k: {}} \
g  \mapsto  (\gamma_k(g), \overline{\gamma_k(\sigma_G(g))}).
\end{gather*}
It can be easily seen that $(\gamma_k\times\overline{\sigma_\mathbb{C}^*}\gamma_k)^*(h_1^G)=\lambda_k^G$.
So $(\lambda_k^G)^2=0$.
\end{proof}

\subsection[The module structure of $KR^*_{(U(n), \sigma_\mathbb{F})}(U(n),\sigma_\mathbb{F})$]{The
module structure of $\boldsymbol{KR^*_{(U(n), \sigma_\mathbb{F})}(U(n),\sigma_\mathbb{F})}$}

\begin{Definition}
Let $\sigma_n$ be (the class of) the standard representation of $U(n)$.
\end{Definition}

\begin{Proposition}
\label{realrepunitary}
$\sigma_n,\wedge^2\sigma_n, \dots, \wedge^n\sigma_n\in RR(U(n), \sigma_\mathbb{R}, \mathbb{R})$, $\wedge^{2i}
\sigma_{2m}
\in RR(U(2m),
\sigma_\mathbb{H}, \mathbb{R})$ and $\wedge^{2i+1}
\sigma_{2m}
\in RH(U(2m),
\sigma_\mathbb{H}, \mathbb{R})$.
Also, both $R(U(n),
\sigma_\mathbb{R}, \mathbb{C})$ and $R(U(2m),
\sigma_\mathbb{H}, \mathbb{C})$ are~$0$.
\end{Proposition}
\begin{proof}
For the involution $\sigma_\mathbb{R}$ and $\wedge^i
\sigma_n$, def\/ine the bilinear form
\begin{gather*}
B_\mathbb{R}: \  \wedge^i\sigma_n\times\sigma_\mathbb{R}^*\wedge^i\sigma_n  \to  \mathbb{C},
\\
\phantom{B_\mathbb{R}: {}} \
(v_1\wedge\dots\wedge v_i, w_1\wedge\dots\wedge w_i)  \mapsto  \det(\langle v_j, \overline{w_k}\rangle).
\end{gather*}
Obviously the form is $U(n)$-invariant, symmetric and non-degenerate.
By Proposition~\ref{formclassification}, each of $\wedge^i
\sigma_n$, $1\leq i\leq n$ is a~Real representation of real type.
Similarly, def\/ine, for the involution $\sigma_\mathbb{H}$ and $\wedge^i
\sigma_{2m}$, a~bilinear form
\begin{gather*}
B_\mathbb{H}: \ \wedge^i\sigma_{2m}\times\sigma_\mathbb{H}^*{\wedge^i\sigma_{2m}}  \to  \mathbb{C},
\\
\phantom{B_\mathbb{H}: {}} \
(v_1\wedge\dots\wedge v_i, w_1\wedge\dots\wedge w_i)  \mapsto  \det(\langle J_mv_j, \overline{w_k}\rangle).
\end{gather*}
It is $U(n)$-invariant because
\begin{gather*}
B_\mathbb{H}(gv, gw)=\det(\langle J_mgv_j, \overline{J_m\overline{g}J_m^{-1}w_k}\rangle)
=\det(\langle J_mgv_j, J_mgJ_m^{-1}\overline{w_k}\rangle)
\\
\phantom{B_\mathbb{H}(gv, gw)}
=\det(\langle v_j, J_m^{-1}\overline{w_k}\rangle)
=\det(\langle J_mv_j, \overline{w_k}\rangle).
\end{gather*}
Moreover
\begin{gather*}
B_\mathbb{H}(v, w)=\det(\langle J_mv_j, \overline{w_k}\rangle)
=\det(\langle -v_j, J_m\overline{w_k}\rangle)
=\det(-\overline{\langle J_m\overline{w_k}, v_j\rangle})
\\
\phantom{B_\mathbb{H}(v, w)}
=\det(-\langle J_mw_k, \overline{v_j}\rangle)
=(-1)^iB_\mathbb{H}(w, v).
\end{gather*}
So by Propositions~\ref{realrepclassification},~\ref{formclassification} and~\ref{quatrepclassification}, $\wedge^i
\sigma_{2m}$ is a~Real representation of real type when $i$ is even and a~Quaternionic representation
of real type when $i$ is odd.
There are no complex representations of complex type because
$\wedge^i\sigma_n\cong \sigma_\mathbb{F}^*\overline{\wedge^i\sigma_n}$ for $\mathbb{F}=\mathbb{R}$ and $\mathbb{H}$.
\end{proof}

\begin{Lemma}
$KR^*_{(U(n),
\sigma_\mathbb{\mathbb{F}})}\!(U(n){,}
\sigma_\mathbb{F})$ is isomorphic to $\Omega^*_{KR^*_{(U(n),
\sigma_\mathbb{F})}(\pt)/ KR^*(\pt)}\!$ as $KR^*_{(U(n),
\sigma_\mathbb{F})}(\pt)$-modules,
where $\mathbb{F}=\mathbb{R}$ or $\mathbb{H}$.
\end{Lemma}
\begin{proof}
Theorem~\ref{equivstrthm} implies that, $KR^*_{(U(n),
\sigma_\mathbb{R})}(U(n),
\sigma_\mathbb{R})\cong RR(U(n),
\sigma_\mathbb{R})\otimes KR^*(U(n),
\sigma_\mathbb{R})$ and $KR^*_{(U(n),
\sigma_\mathbb{H})}(U(n),
\sigma_\mathbb{H})\cong (RR(U(n),
\sigma_\mathbb{H}, \mathbb{R})\oplus RH(U(n),
\sigma_\mathbb{H}, \mathbb{R}))\otimes KR^*(U(n),
\sigma_\mathbb{H})$.
Moreover, by Theorem~\ref{spkrtheorymodofg} and Proposition~\ref{realrepunitary},
\begin{gather*}
KR^*(U(n),
\sigma_\mathbb{R})\cong \bigwedge\nolimits_{KR^*(\pt)}(\delta_\mathbb{R}(\sigma_n), \dots, \delta_\mathbb{R}
(\wedge^n
\sigma_n))
\end{gather*}
and
\begin{gather*}
KR^*(U(2m),
\sigma_\mathbb{H})\cong\bigwedge\nolimits_{KR^*(\pt)}\big(\delta_\mathbb{H}(\sigma_{2m}), \delta_\mathbb{R}(\wedge^2
\sigma_{2m}), \dots, \delta_\mathbb{R}(\wedge^{2m}\sigma_{2m})\big).
\end{gather*}
Putting all these together and applying Theorem~\ref{KR_G}, we get the desired conclusion.
\end{proof}
\begin{Remark}
\label{ungradediso}
As ungraded $KR^*(\pt)$-modules, both
\[
KR_{(U(2m),
\sigma_\mathbb{R})}^*(U(2m),
\sigma_\mathbb{R}) \qquad \text{and}\qquad KR^*_{(U(2m),
\sigma_\mathbb{H})}(U(2m),
\sigma_\mathbb{H})
\]
 are isomorphic to $K^*_{U(2m)}(U(2m))\otimes KR^*(\pt)$.
\end{Remark}

\subsection{Injectivity results}
This step involves proving that the restriction map $p_G^*$ to the equivariant $KR$-theory of the maximal torus and the
map $q_G^*$ induced by the Weyl covering map are injective.
\begin{Lemma}
\label{ringinj}
Let~$G$ be a~compact Lie group and $X$ a~$G$-space.
Let $i_1^*: K_G^*(X)\to K^*_T(X)$ be the map which restricts the~$G$-action to $T$-action.
Then
\begin{gather*}
i_1^*\otimes\Id_R: \  K^*_G(X)\otimes R\to K_T^*(X)\otimes R
\end{gather*}
is injective for any ring $R$.
\end{Lemma}
\begin{proof}
By~\cite[Proposition 4.9]{At4}, $i_1^*$ is split injective.
So is $i_1^*\otimes\Id_R$ for any ring $R$.
\end{proof}
\begin{Lemma}
\label{decomposition}
Let $i_2^*: K_T^*(G)\to K_T^*(T)\cong R(T)\otimes K^*(T)$ be the map induced by the inclusion $T\hookrightarrow G$.
Then
\begin{gather*}
i_2^*\delta_T(\rho)=
\sum\limits_{j=1}^{\dim\rho}e^{\tau_j}\otimes\delta(\tau_j)\in K_T^{-1}(T),
\end{gather*}
where $\tau_j$ are the weights of $\rho$.
\end{Lemma}
\begin{proof}
Let~$V$ be the vector space underlying the representation $\rho$.
$\delta_T(\rho)$ is represented by the complex of $T$-equivariant vector bundles
\begin{gather*}
0\longrightarrow G\times\mathbb{R}\times V \longrightarrow G\times\mathbb{R}\times V\longrightarrow 0,
\\
(g, t, v)  \mapsto  (g, t, -t\rho(g)v) \qquad\text{if} \ \ t\geq 0,
\\
(g, t, v)  \mapsto  (g, t, tv) \qquad\text{if} \ \ t\leq 0
\end{gather*}
which, on restricting to
\begin{gather*}
0\longrightarrow T\times\mathbb{R}\times V
\stackrel{}
{\longrightarrow}T\times\mathbb{R}\times V\longrightarrow 0
\end{gather*}
is decomposed into a~direct sum of complexes of 1-dimensional $T$-equivariant vector bundles, each of which corresponds
to a~weight of $\rho$.
\end{proof}

\begin{Lemma}
\label{Atiyahindex}
Let~$G$ be a~simply-connected, connected compact Lie group and $\rho_1, \dots, \rho_l$ be its fundamental
representations.
Then
\begin{gather}
\label{jacobian}
i_2^*\left(\prod\limits_{i=1}^l\delta_T(\rho_i)\right)=d_G\otimes\prod\limits_{i=1}^l\delta(\varpi_i),
\end{gather}
where $\varpi_i$ the $i$-th fundamental weight and $d_G=
\sum\limits_{w\in W}
\operatorname{sgn}(w)e^{w\cdot
\sum\limits_{i=1}^l\varpi_i}\in R(T)$ is the Weyl deno\-mi\-nator.
\end{Lemma}
\begin{proof}
Equation~\eqref{jacobian} follows from Lemma~\ref{decomposition} and Lemma 3 of~\cite{At2}.
\end{proof}
\begin{Lemma}
\label{wedgeinj}
Let $M$ be an $R$-module freely generated by $m_1, \dots, m_l$, and $N$ an $R$-module.
If $f: M\to N$ is an $R$-module homomorphism, and $rf(m_1)\wedge\dots\wedge f(m_l)\in \bigwedge_R^l N$ is nonzero for
all $r\in R
\setminus\{0\}$, then
\begin{gather*}
\bigwedge\nolimits^* f: \ \bigwedge\nolimits_R^*M\to \bigwedge\nolimits_R^*N
\end{gather*}
is injective.
\end{Lemma}
\begin{proof}
It suf\/f\/ices to show that $\bigwedge^k f$ is injective for $1\leq k\leq l$.
Suppose $I
\subseteq\{1, \dots, l\}$, $|I|=k$, $m_I:=\bigwedge_{i\in I}m_i$ and $f(m_I):=\bigwedge_{i\in I}f(m_i)$.
If $\sum\limits_{|I|=k}
r_Im_I\in\text{ker}(\bigwedge^k f)$, then for any $J$ with $|J|=k$,
\begin{gather*}
0=
\sum\limits_{|I|=k}
r_If(m_I)\wedge f(m_{J^c})=r_Jf(m_1)\wedge\dots\wedge f(m_l).
\end{gather*}
Hence $\sum\limits_{|I|=k}
r_Im_I=0$ and the conclusion follows.
\end{proof}
\begin{Lemma}
\label{twotorsioninj}
Let~$G$ be a~simply-connected, connected and compact Lie group.
Then the map
\begin{gather*}
i_2^*\otimes\Id_R: \ K_T^*(G)\otimes R\to K_T^*(T)\otimes R
\end{gather*}
is injective for any ring $R$.
\end{Lemma}
\begin{proof}
Note that
\begin{gather*}
K_G^*(G)\otimes_{R(G)}R(T)  \to  K_T^*(G),
\\
a\otimes\rho  \mapsto  i_1^*(a)\cdot\rho.
\end{gather*}
is an $R(T)$-algebra isomorphism (cf.~\cite[Theorem~4.4]{HLS}).
Using Theorem~\ref{eqderivation}, we have that $K_T^*(G)$ is isomorphic, as an $R(T)$-algebra, to $\bigwedge_{R(T)}^*M$,
where $M$ is the $R(T)$-module freely generated by $\delta_T(\rho_1), \dots, \delta_T(\rho_l)$.
We also observe that $K_T^*(T)$ is isomorphic, as an $R(T)$-algebra, to $\bigwedge_{R(T)}^*N$, where $N$ is the
$R(T)$-module freely generated by $\delta(\varpi_1), \dots, \delta(\varpi_l)$.
Note that the hypotheses of Lemma~\ref{wedgeinj} are satisf\/ied by $f=i_2^*\otimes\Id_{\mathbb{Z}_m}$ for any
$m\geq 2$, as
$r\prod\limits_{i=1}^li_2^*\delta_T(\rho_i)=ri_2^*\prod\limits_{i=1}^l\delta_T(\rho_i)=rd_G\otimes\prod\limits_{i=1}^l\delta(\varpi_i)$
(by Lemma~\ref{Atiyahindex}) is indeed nonzero for any nonzero $r$ in $\mathbb{Z}_m$ (the coef\/f\/icients of $d_G$ are
either 1 or $-1$, so after reduction mod $m$ $rd_G$ is still nonzero).
Now that $i_2^*\otimes\mathbb{Z}_m$ is injective, so is $i_2^*\otimes\Id_R$ for any ring $R$.
\end{proof}

\begin{Proposition}\label{pginjective}{\samepage\quad
\begin{enumerate}\itemsep=0pt
\item[$1.$]
Let~$G$ be a~compact, connected and simply-connected Real Lie group such that $RR(G, \mathbb{C})=RR(G, \mathbb{H})=0$
and there exists a~maximal torus $T$ on which the involution acts by inversion.
Then the restriction map
$KR_{G}^*(G)\to KR_T^*(T)$
is injective.
\item[$2.$]
The map $p_G^*: KR_{(U(n),
\sigma_\mathbb{R})}^*(U(n),
\sigma_\mathbb{R})\to KR_{(T,
\sigma_\mathbb{R})}^*(T,
\sigma_\mathbb{R})$ is injective.
\end{enumerate}}
\end{Proposition}

\begin{proof}
By Corollary~\ref{equivstrthmg}, $KR^*_G(G)\cong K^*_G(G)\otimes KR^*(\pt)$ and $KR_T^*(T)\cong K^*_T(T)\otimes
KR^*(\pt)$ as $KR^*(\pt)$-modules.
Using this identif\/ication, we can as well identify the restriction map with $i^*\otimes\Id_{KR^*(\pt)}$,
where $i^*:=i_2^*\circ i_1^*$.
Part~(1) then follows from Lemmas~\ref{ringinj} and~\ref{twotorsioninj}.
Part~(2) is immediate once we apply Lemma~\ref{ringinj} and note that the proofs of Lemmas~\ref{Atiyahindex}
and~\ref{twotorsioninj} can be adapted to the case $G=U(n)$ by letting $\sigma_n, \dots, \wedge^n\sigma_n$
play the role of the fundamental representations and their highest weights, the
fundamental weights.
\end{proof}

\begin{Lemma}
\label{quatring}
$KR^*_{(U(2m),
\sigma_\mathbb{H})}(U(2m)/T,
\sigma_\mathbb{H})\cong \mathbb{Z}[e_1^\mathbb{H}, \dots, e_{2m}^\mathbb{H}, (e_1^\mathbb{H}e_2^\mathbb{H}\cdots
e_{2m}^\mathbb{H})^{-1}]\otimes KR^*(\pt)$ as rings, where $e_i^\mathbb{H}$ lives in the degree $-4$ piece.
\end{Lemma}
\begin{proof}
It is known that $K^*(U(2m)/T)\cong\mathbb{Z}[\alpha_1, \dots, \alpha_{2m}]/(s_i-\binom{2m}{i}|1\leq i\leq 2m)$, where
$\alpha_i=[U(2m)\times_T\mathbb{C}_{e_i}]$ and $s_i$ is the $i$-th elementary symmetric polynomial (cf.~\cite[Proposition 2.7.13]{At}).
The induced map $\overline{\sigma_\mathbb{H}^*}$ acts as identity on $K^*(U(2m)/T)$.
The involution $\sigma_\mathbb{H}$ on the base lifts to a~Quaternionic structure
on the associated complex line bundle, so there \mbox{exist}
$\alpha_1^\mathbb{H}, \dots, \alpha_{2m}^\mathbb{H}\in KR^{-4}(U(2m)/T)$, such that their complexif\/ications are
$\beta^2\alpha_1, \dots, \beta^2\alpha_{2m}\in K^*(U(2m)/T)$.
By Theorem~\ref{strthm}, $KR^*(U(2m)/T,
\sigma_\mathbb{H})$ is a~$KR^*(\pt)$-module generated by polynomials in $\alpha_1^\mathbb{H}, \dots, \alpha_{2m}^\mathbb{H}\in
KR^{-4}(U(2m)/T,
\sigma_\mathbb{H})$.
In fact it is not hard to see that $KR^*(U(2m)/T$,
$\sigma_\mathbb{H})$ is isomorphic to
\begin{gather*}
\mathbb{Z}[\alpha_1^\mathbb{H}, \dots, \alpha_{2m}^\mathbb{H}]\otimes
KR^*(\pt)\left/\left.\left(s_{2k}-\binom{2m}{2k}, s_{2k-1}-\frac{1}{2}\mu\binom{2m}{2k-1}\right|1\leq k\leq
m\right)\right..
\end{gather*}
Also obvious is that each of $\alpha_i^\mathbb{H}$ has an equivariant lift $e_i^\mathbb{H}\in KR^{-4}_{(U(2m),
\sigma_\mathbb{H})}(U(2m)/T,\sigma_\mathbb{H})$.

Now that all the hypotheses in Theorem~\ref{equivstrthm} are satisf\/ied, we can apply it, together with the fact that
$R(U(2m),
\sigma_\mathbb{H}, \mathbb{C})=0$ (cf.\ Proposition~\ref{realrepunitary}) to see that $KR^*_{(U(2m),
\sigma_\mathbb{H})}(U(2m)/T,
\sigma_\mathbb{H})$ is isomorphic to
\begin{gather*}
(RR(U(2m),
\sigma_\mathbb{H}, \mathbb{R})\oplus RH(U(2m),
\sigma_\mathbb{H}, \mathbb{R}))\otimes KR^*(U(2m)/T,
\sigma_\mathbb{H})
\end{gather*}
as $RR(U(2m),
\sigma_\mathbb{H}, \mathbb{R})\oplus RH(U(2m),
\sigma_\mathbb{H}, \mathbb{R})$-modules (actually as rings).
Noting that
\begin{gather*}
RR(U(2m),\sigma_\mathbb{H}, \mathbb{R})\oplus RH(U(2m),
\sigma_\mathbb{H}, \mathbb{R})\cong\mathbb{Z}\big[s_1, \dots, s_{2m}, s_{2m}^{-1}\big]
\end{gather*}
we establish the Lemma.
\end{proof}

\begin{Proposition}
$KR_{(U(2m),
\sigma_\mathbb{H})}^*(U(2m)/T\times T,
\sigma_\mathbb{H}
\times
\sigma_\mathbb{R})$ is isomorphic to
\begin{gather*}
\mathbb{Z}\big[e_1^\mathbb{H}, \dots, e_{2m}^\mathbb{H},
\big(e_1^\mathbb{H}\cdots e_{2m}^\mathbb{H}\big)^{-1}\big]
\otimes KR^*(T,\sigma_\mathbb{R})
\end{gather*}
as graded rings.
\end{Proposition}
\begin{proof}
First, by~\cite[Theorem~1]{P}, Proposition~\ref{KR_G} and Lemma~\ref{quatring}, $KR^*_{(U(2m),
\sigma_\mathbb{H})}(U(2m)/T,
\sigma_\mathbb{H})$ is a~free $KR_{(U(2m),
\sigma_\mathbb{H})}^*(\pt)$-module.
The same is also true of $KR_{(U(2m),
\sigma_\mathbb{H})}^*(T,
\sigma_\mathbb{R})$ since, by Theo\-rem~\ref{equivstrthm}, it is isomorphic to $RR(U(2m),
\sigma_\mathbb{H})\otimes KR^*(T,
\sigma_\mathbb{R})$, which in turn is isomorphic to $KR^*_{(U(2m),
\sigma_\mathbb{H})}(\pt)\otimes K^*(T)$.
The proposition
follows from a~version of K\"unneth formula for equivariant $KR$-theory.
\end{proof}
\begin{Remark}
\label{ungradedisoh}
$KR^*_{(U(2m),
\sigma_\mathbb{H})}(U(2m)/T\times T,
\sigma_\mathbb{H}
\times
\sigma_\mathbb{R})$ is isomorphic to $K^*_T(T)\otimes KR^*(\pt)$ and $KR^*_{(T,
\sigma_\mathbb{R})}(T,
\sigma_\mathbb{R})$ as ungraded $KR^*(\pt)$-modules.
\end{Remark}
\begin{Proposition}
\label{qginjective}
The map $q_G^*$ is injective.
\end{Proposition}
\begin{proof}
By Remarks~\ref{ungradediso} and~\ref{ungradedisoh}, $KR^*_{(U(2m),
\sigma_\mathbb{H})}(U(2m),
\sigma_\mathbb{H})$ and $KR^*_{(U(2m),
\sigma_\mathbb{H})}(U(2m)/T\times T,
\sigma_\mathbb{H}
\times
\sigma_\mathbb{R})$ are isomorphic to $KR^*_{(U(2m),
\sigma_\mathbb{R})}(U(2m),
\sigma_\mathbb{R})$ and $KR^*_{(T,
\sigma_\mathbb{R})}(T,
\sigma_\mathbb{R})$ respectively, as ungraded $KR^*(\pt)$-modules.
It is not hard to see that $q_G^*$ can be identif\/ied with $p_G^*$ under these isomorphisms.
Now the result follows from Proposition~\ref{pginjective}.
\end{proof}
\subsection{Squares of algebra generators of real and quaternionic types}
\begin{Lemma}
\label{eqkrthytorus}
$KR^*(T,
\sigma_\mathbb{R})$ is isomorphic to the exterior algebra over $KR^*(\pt)$ generated by $\delta_\mathbb{R}(e_1), \dots,
\delta_\mathbb{R}(e_n)$, as $KR^*_{(T,
\sigma_\mathbb{R})}(\pt)$-modules.
Here $e_i$ is the $1$-dimensional complex representation with weight being the $i$-th standard basis vector of the weight
lattice.
Moreover, $\delta_\mathbb{R}(e_i)^2=\eta\delta_\mathbb{R}(e_i)$.
\end{Lemma}

\begin{proof}
Since $R(T,
\sigma_\mathbb{R}, \mathbb{C})=0$, the module structure follows from Theorem~\ref{spkrtheorymodofg}.
For the second part of the Lemma, see the appendix of~\cite{Se}.
\end{proof}
\begin{Proposition}
\label{equivkrthysquaresigma}
In $KR^*_{(U(n),
\sigma_\mathbb{R})}(U(n),
\sigma_\mathbb{R})$
\begin{gather}
\label{realsigmasquare}
\delta_\mathbb{R}^G(\sigma_n)^2
=\eta\big(\sigma_n\cdot\delta_\mathbb{R}^G(\sigma_n)-\delta_\mathbb{R}^G(\wedge^2\sigma_n)\big).
\end{gather}
In $KR^*_{(U(2m),
\sigma_\mathbb{H})}(U(2m),
\sigma_\mathbb{H})$,
\begin{gather}
\label{quatsigmasquare}
\delta_\mathbb{H}^G(\sigma_{2m})^2=\eta\big(\sigma_{2m}
\cdot\delta_\mathbb{H}^G(\sigma_{2m})-\delta_\mathbb{R}^G(\wedge^2\sigma_{2m})\big).
\end{gather}
\end{Proposition}
\begin{proof}
Now that we have shown that $p^*_G$ and $q^*_G$ are injective by Propositions~\ref{pginjective} and~\ref{qginjective},
we can compute $\delta_\mathbb{R}^G(\sigma_n)^2$ and $\delta_\mathbb{H}^G(\sigma_{2m})^2$
by passing the computation through $p^*_G$ and $q_G^*$ to their images.
We prove the case $\mathbb{F}=\mathbb{R}$.
The proof of the case $\mathbb{F}=\mathbb{H}$ is similar so we leave it to the reader.
Note that
\begin{gather*}
p^*_G\big(\delta_\mathbb{R}^G(\sigma_n)^2\big)=p^*_G(\delta_\mathbb{R}^G(\sigma_n))^2
=\left(\sum\limits_{i=1}^n e_i\otimes\delta_\mathbb{R}(e_i)\right)^2
\\
\phantom{p^*_G\big(\delta_\mathbb{R}^G(\sigma_n)^2\big)}=
\sum\limits_{i=1}^n e_i^2\otimes\eta\delta_\mathbb{R}(e_i)+
\sum\limits_{i\neq j}
e_ie_j\otimes\delta_\mathbb{R}(e_i)\delta_\mathbb{R}(e_j)
\\
\phantom{p^*_G\big(\delta_\mathbb{R}^G(\sigma_n)^2\big)} \overset{\text{Lemma~\ref{eqkrthytorus}}}{=}
\sum\limits_{i=1}^ne_i^2\otimes\eta \delta_\mathbb{R}(e_i)+
\sum\limits_{i<j}
e_i e_j\otimes\big(\delta_\mathbb{R}(e_i)\delta_\mathbb{R}(e_j)+\delta_\mathbb{R}(e_j)\delta_\mathbb{R}(e_i)\big)
\\
\phantom{p^*_G\big(\delta_\mathbb{R}^G(\sigma_n)^2\big)}=
\sum\limits_{i=1}^n e_i^2\otimes\eta\delta_\mathbb{R}(e_i).
\end{gather*}
On the other hand,
\begin{gather*}
p^*(\sigma_n\cdot\delta_\mathbb{R}^G(\sigma_n))=\left(\sum\limits_{i=1}^ne_i\otimes 1\right)
\left(\sum\limits_{i=1}^n e_i\otimes\delta_\mathbb{R}(e_i)\right)
=
\sum\limits_{1\leq i, j\leq n}
e_ie_j\otimes\delta_\mathbb{R}(e_i)
\\
\phantom{p^*(\sigma_n\cdot\delta_\mathbb{R}^G(\sigma_n))}=
\sum\limits_{i=1}^n e_i^2\otimes \delta_\mathbb{R}(e_i)+
\sum\limits_{i\neq j}
e_j e_i\otimes\delta_\mathbb{R}(e_i)
\\
\phantom{p^*(\sigma_n\cdot\delta_\mathbb{R}^G(\sigma_n))}=
\sum\limits_{i=1}^n e_i^2\otimes\delta_\mathbb{R}(e_i)+
\sum\limits_{1\leq i<j\leq n}
e_j e_i\otimes(\delta_\mathbb{R}(e_i)+\delta_\mathbb{R}(e_j))
\\
\phantom{p^*(\sigma_n\cdot\delta_\mathbb{R}^G(\sigma_n))}=
\sum\limits_{i=1}^n e_i^2\otimes\delta_\mathbb{R}(e_i)+p^*_G\big(\delta_\mathbb{R}^G(\wedge^2 \sigma_n)\big).
\end{gather*}
From the above equations we obtain
\begin{gather*}
p^*_G\big(\delta_\mathbb{R}^G(\sigma_n)^2\big)
=\eta p^*_G(\sigma_n\cdot\delta_\mathbb{R}^G(\sigma_n)-\delta_\mathbb{R}^G(\wedge^2\sigma_n)).
\end{gather*}
By Proposition~\ref{pginjective},
\begin{gather*}
\delta_\mathbb{R}^G(\sigma_n)^2
=\eta\big(\sigma_n\cdot\delta_\mathbb{R}^G(\sigma_n)-\delta_\mathbb{R}^G(\wedge^2\sigma_n)\big). \tag*{\qed}
\end{gather*}
\renewcommand{\qed}{}
\end{proof}
\begin{Theorem}
\label{equivkrthysquare}
Let~$G$ be a~Real compact Lie group.
Then
\begin{gather}
\label{realsquare}
\delta_\mathbb{R}^G(\varphi_i)^2
=\eta\big(\varphi_i\cdot\delta_\mathbb{R}^G(\varphi_i)-\delta_\mathbb{R}^G(\wedge^2\varphi_i)\big),
\\
\label{quatsquare}
\delta_\mathbb{H}^G(\theta_j)^2
=\eta\big(\theta_j\cdot\delta_\mathbb{H}^G(\theta_j)-\delta_\mathbb{R}^G(\wedge^2\theta_j)\big).
\end{gather}
\end{Theorem}
\begin{proof}
The induced map $\varphi_i^*: KR^*_{(U(n),
\sigma_\mathbb{R})}(U(n),
\sigma_\mathbb{R})\to KR^*_{(G,
\sigma_G)}
(G,
\sigma_G)$ sends $\sigma_n$ to $\varphi_i$, and $\delta_\mathbb{R}^G(\sigma_n)$ to $\delta_\mathbb{R}^G(\varphi_i)$.
Likewise, the induced map $\theta_j^*: KR^*_{(U(2m),
\sigma_\mathbb{H})}(U(2m),
\sigma_\mathbb{H})$ sends $\sigma_{2m}$ to $\theta_j$, and $\delta_\mathbb{H}^G(\sigma_{2m})$ to $\delta_\mathbb{H}^G(\theta_j)$.
The result now follows from Proposition~\ref{equivkrthysquaresigma}.
\end{proof}

To further express $\delta_\mathbb{R}^G(\wedge^2\varphi_i)$ and $\delta_\mathbb{R}^G(\wedge^2\theta_j)$ in terms of the
module generators associated with the fundamental representations, we may use the following derivation property of
$\delta_\mathbb{R}^G$ and $\delta_\mathbb{H}^G$.
\begin{Proposition}
\label{eqderivative}
$\delta_\mathbb{R}^G\oplus\delta_\mathbb{H}^G$ is a~derivation of the graded ring $RR(G)\oplus RH(G)$ $($with $RR(G)$ of
degree $0$ and $RH(G)$ of degree $-4)$ taking values in the graded module $KR^{-1}_G(G)\oplus KR^{-5}_G(G)$.
\end{Proposition}
\begin{proof}
We refer the reader to the proof of Proposition~3.1 of~\cite{BZ}
with the def\/inition of $\delta_G(\rho)$ given there
(which is incorrect) replaced by the one in Def\/inition~\ref{defdeltag}.
One just need to simply check that the homotopy $\rho_s$ in the proof for $t\geq 0$ intertwines with both $\sigma_\mathbb{R}$ and $\sigma_\mathbb{H}$.
\end{proof}
\begin{Corollary}
\label{krringunitary}
In $KR^*_{(U(n),
\sigma_\mathbb{R})}(U(n),
\sigma_\mathbb{R})$,
\begin{gather*}
\delta_\mathbb{R}^G\big({\wedge}^k
\sigma_n\big)^2=\eta\sum\limits_{i=1}^{2k}\wedge^{2k-i}
\sigma_n\cdot\delta_\mathbb{R}^G\big({\wedge}^i
\sigma_n\big).
\end{gather*}
In $KR^*_{(U(2m),
\sigma_\mathbb{H})}(U(2m),
\sigma_\mathbb{H})$,
\begin{gather*}
\delta_\mathbb{H}^G\big({\wedge}^{2k-1}\sigma_{2m}\big)^2=\eta\sum\limits_{j=1}^{2k-1}\left(\wedge^{4k-2j-1}
\sigma_{2m}\cdot\delta_\mathbb{H}^G\big({\wedge}^{2j-1}\sigma_{2m}\big)+\wedge^{4k-2j-2}\sigma_{2m}
\cdot\delta_\mathbb{R}^G\big({\wedge}^{2j}\sigma_{2m}\big)\right),
\\
\delta_\mathbb{R}^G\big({\wedge}^{2k}\sigma_{2m}\big)^2=\eta\sum\limits_{j=1}^{2k}\left(\wedge^{4k-2j+1}
\sigma_{2m}\cdot\delta_\mathbb{H}^G\big({\wedge}^{2j-1}\sigma_{2m}\big)+\wedge^{4k-2j}\sigma_{2m}
\cdot\delta_\mathbb{R}^G\big({\wedge}^{2j}\sigma_{2m}\big)\right).
\end{gather*}
\end{Corollary}
\begin{proof}
By the def\/inition,
\begin{gather*}
\big({\wedge}^k
\sigma_n\big)^*\big(\delta_\mathbb{R}^G\big({\wedge}^2
\sigma_{\left(
\begin{array}{@{}c@{}}
n
\\
k
\end{array}\right)
}\big)\big)=\delta_\mathbb{R}^G\big({\wedge}^2\big({\wedge}^k
\sigma_n\big)\big).
\end{gather*}
By Exercise 15.32 of~\cite{FH},
\begin{gather*}
\wedge^2(\wedge^k
\sigma_n)=\bigoplus_i\Gamma_{\varpi_{k-2i+1}
+\varpi_{k+2i-1}}.
\end{gather*}
By Giambelli's formula,
\begin{gather*}
\Gamma_{\varpi_{k-2i+1}+\varpi_{k+2i-1}}=\wedge^{k+2i-1}
\sigma_n\cdot\wedge^{k-2i+1}\sigma_n-\wedge^{k+2i}\sigma_n\cdot\wedge^{k-2i}\sigma_n.
\end{gather*}
By Proposition~\ref{eqderivative}
\begin{gather*}
\delta_\mathbb{R}^G(\Gamma_{\varpi_{k-2i+1}+\varpi_{k+2i-1}})=\wedge^{k+2i-1}
\sigma_n\cdot\delta_\mathbb{R}^G\big({\wedge}^{k-2i+1}\sigma_n\big)+\wedge^{k-2i+1}\sigma_n\cdot\delta_\mathbb{R}^G\big({\wedge}^{k+2i-1}\sigma_n\big)
\\
\hphantom{\delta_\mathbb{R}^G(\Gamma_{\varpi_{k-2i+1}+\varpi_{k+2i-1}})=}
{}-\wedge^{k+2i}
\sigma_n\cdot\delta_\mathbb{R}^G\big({\wedge}^{k-2i}\sigma_n\big)-\wedge^{k-2i}
\sigma_n\cdot\delta_\mathbb{R}^G\big({\wedge}^{k+2i}\sigma_n\big).
\end{gather*}
Now the f\/irst equation is immediate.
The second and third equations can be derived similarly.\looseness=1
\end{proof}

Putting together the previous results yields the following full description of the ring structure of $KR_G^*(G)$.
\begin{Theorem}
\label{mainthm2}
Let~$G$ be a~simply-connected, connected and compact Real Lie group.
Viewing~$G$ as a~Real~$G$-space with adjoint action, we have
\begin{enumerate}\itemsep=0pt
\item[$1.$] {\rm (}Corollary~{\rm \ref{equivstrthmg})} The map
\begin{gather*}
f: \ (RR(G, \mathbb{R})\oplus RH(G, \mathbb{R}))\otimes KR^*(G)\oplus r(R(G, \mathbb{C})\otimes K^*(G))\to KR_G^*(G)
\\
\phantom{f: {}} \
\rho_1\otimes a_1\oplus r(\rho_2\otimes a_2)\mapsto \rho_1\cdot (a_1)_{G}\oplus r(\rho_2\cdot (a_2)_G)
\end{gather*}
is a~group isomorphism.
In particular, if $R(G, \mathbb{C})=0$, then $f$ is an isomorphism of $KR^*_{G}(\pt)$-modules.
\item[$2.$] {\rm (}Corollary~{\rm \ref{eqmodgenerators})}
$KR_{G}(G)$ is generated by $\delta_\mathbb{R}^G(\varphi_1), {\dots}, \delta_\mathbb{R}^G(\varphi_r)$,
$\delta_\mathbb{H}^G(\theta_1), {\dots}, \delta_\mathbb{H}^G(\theta_s)$, $\lambda_1^G, {\dots}, \lambda_t^G\!$ and
$r^G_{\rho, i, \varepsilon_1, \dots, \varepsilon_t, \nu_1, \dots, \nu_t}$ as an algebra over $KR_{G}^*(\pt)$.
Moreover,
\begin{gather}
(a)\quad{\rm (}Proposition~{\rm \ref{cplxsquare})} \quad \big(\lambda_k^G\big)^2=0\text{ for all }1\leq k\leq t,\nonumber
\\ \label{cplxrel1}
(b)\quad
(r^G_{\rho, i, \varepsilon_1, \dots, \varepsilon_t, \nu_1, \dots, \nu_t})^2
\\
\phantom{(b)\quad}
=
\begin{cases}
\eta^2(\rho\cdot\overline{\sigma_G^*}
\rho)(\lambda^G_1)^{\omega_1}\cdots(\lambda_t^G)^{\omega_t}
&
\text{if }r^G_{\rho, i, \varepsilon_1, \dots, \varepsilon_t, \nu_1, \dots, \nu_t}
\text{ is of degree} \ -1\text{ or }-5,
\\
\pm\mu(\rho\cdot\overline{\sigma_G^*}
\rho)(\lambda^G_1)^{\omega_1}\cdots(\lambda_t^G)^{\omega_t}
&
\text{if }r^G_{\rho, i, \varepsilon_1, \dots, \varepsilon_t, \nu_1, \dots, \nu_t}
\text{ is of degree} \ -2\text{ or }-6,
\\
0
&\text{otherwise}.
\end{cases}
\nonumber
\end{gather}
The sign can be determined using formulae in $(2)$ of Proposition~{\rm \ref{krprelim}}.
\begin{gather}
\label{cplxrel2}
(c)\quad
r^G_{\rho, i, \varepsilon_1, \dots, \varepsilon_t, \nu_1, \dots, \nu_t}\eta=0,\text{ and }
r^G_{\rho, i, \varepsilon_1,\dots, \varepsilon_t, \nu_1, \dots, \nu_t}\mu
=2r^G_{\rho, i+2, \varepsilon_1, \dots, \varepsilon_t, \nu_1, \dots,\nu_t},
\\
(d)\quad
{\rm (}Proposition~{\rm \ref{equivkrthysquare})}
\quad
\delta_\mathbb{R}^G(\varphi_i)^2
=\eta\big(\varphi_i\cdot\delta_\mathbb{R}^G(\varphi_i)-\delta_\mathbb{R}^G(\wedge^2\varphi_i)\big),
\nonumber
\\
\phantom{(d)\quad}
\delta_\mathbb{H}^G(\theta_j)^2
=\eta\big(\theta_j\cdot\delta_\mathbb{H}^G(\theta_j)-\delta_\mathbb{R}^G(\wedge^2\theta_j)\big).
\nonumber
\end{gather}
One can express $\delta_\mathbb{R}^G(\wedge^2\varphi_i)$ and $\delta_\mathbb{H}^G(\wedge^2\theta_j)$ in terms of the
algebra generators using the derivation property of $\delta_\mathbb{R}^G$ and $\delta_\mathbb{H}^G$
$($cf.\ Proposition~{\rm \ref{eqderivative})}.
\end{enumerate}
\end{Theorem}
\begin{proof}
Only~\eqref{cplxrel1} and~\eqref{cplxrel2} need explanation, but they are just equivariant analogues of
Corollary~\ref{krglooseend} and follow from (2) of Proposition~\ref{krprelim} and Remark~\ref{cplxeqsquare}.
\end{proof}
\begin{Remark}{\samepage\quad
\begin{enumerate}\itemsep=0pt
\item
$KR^*_G(\pt)_2$, which is the ring obtained by inverting the prime 2 in $KR^*_G(\pt)$, is isomorphic to\looseness=-1
\begin{gather*}
(RR(G,\sigma_G, \mathbb{R})\oplus RH(G,
\sigma_G, \mathbb{R}))\otimes\mathbb{Z}\left[\tfrac{1}{2}, \beta^2\right]/\big(\big(\beta^2\big)^2-1\big)
\\
\qquad{}
\oplus r\big(R\big(G,\sigma_G, \mathbb{C}\big)
\otimes\mathbb{Z}\left[\tfrac{1}{2}, \beta\right]/(\beta^4-1)\big).
\end{gather*}
\item
If $R(G, \mathbb{C})=0$, then $KR^*_G(G)_2$, which is the ring obtained by inverting the prime 2 in $KR^*_G(G)$, is
isomorphic to, as $KR^*_G(\pt)_2$-algebra,
\begin{gather*}
\bigwedge\nolimits_{KR^*_G(\pt)_2}\left(\delta_\mathbb{R}^G(\varphi_1), \dots, \delta_\mathbb{R}^G(\varphi_r),
\delta_\mathbb{H}^G(\theta_1), \dots, \delta_\mathbb{H}^G(\theta_s)\right).
\end{gather*}
\end{enumerate}}
\end{Remark}

\section{Applications and examples}\label{Section5}

Applying the forgetful map $KR_G^*(G)\to KR^*(G)$ to Proposition~\ref{eqderivative} and Theorem~\ref{mainthm2}, we solve
the problem of f\/inding a~description of the ring structure of $KR^*(G)$ which was left open by Seymour in~\cite{Se}.
\begin{Theorem}
Let~$G$ be a~simply-connected, connected and compact Real Lie group.
Then
\begin{gather*}
\delta_\mathbb{R}(\varphi_i)^2
=\eta\big(\dim(\varphi_i)\cdot\delta_\mathbb{R}(\varphi_i)-\delta_\mathbb{R}(\wedge^2\varphi_i)\big),
\qquad
\delta_\mathbb{H}(\theta_j)^2=\eta\delta_\mathbb{R}(\wedge^2\theta_j).
\end{gather*}
One can express $\delta_\mathbb{R}(\wedge^2\varphi_i)$ and $\delta_\mathbb{H}(\wedge^2\theta_j)$ in terms of the
generators in Proposition~{\rm \ref{krtheorymodofg}} using the derivation property of $\delta_\mathbb{R}$ and
$\delta_\mathbb{H}$ got by applying the forgetful map to Proposition~{\rm \ref{eqderivative}}.
The above equations, together with Theorem~{\rm \ref{krtheorymodofg}} and Corollary~{\rm \ref{krglooseend}}, constitute a~complete
description of the ring structure of~$KR^*(G)$.
\end{Theorem}
\begin{Remark}
Seymour's conjecture concerning $\delta_\mathbb{R}(\sigma_n)^2$ is true.
However, his conjecture that if $x\in KR^{-5}
(X)$, then $x^2=0$ is false, as evidenced by the ring structure of $KR^*(U(2m),
\sigma_{\mathbb{H}})$.
\end{Remark}
\begin{Example}
Let~$G$ be a~simply-connected, connected and compact Real Lie group with no fundamental representations of complex type.
Equip~$G$ with both the trivial~$G$-action and the adjoint action.
Both $K_{G_{\text{triv}}}^*(G)$ and $K_{G_{\text{Ad}}}^*(G)$ are isomorphic to $\Omega^*_{R(G)/\mathbb{Z}}$ as rings.
On the other hand, though both $KR^*_{G_{\text{triv}}}(G)$ and $KR^*_{G_{\text{Ad}}}(G)$ are isomorphic to
$\Omega_{KR^*_G(\pt)/KR^*(\pt)}^*$ as $KR^*_{G}(\pt)$-modules, they are not isomorphic as rings, as
one can tell from the squares of the generators of both rings.
For instance, in $KR_{G_{\text{triv}}}^*(G)$,
\begin{gather*}
\delta^G_\mathbb{R}(\varphi_i)^2
=\eta\left(\dim(\varphi_i)\delta_\mathbb{R}^G(\varphi_i)-\delta_\mathbb{R}^G(\wedge^2\varphi_i)\right)
\end{gather*}
whereas in $KR^*_{G_{\text{Ad}}}(G)$,
\begin{gather*}
\delta_\mathbb{R}^G(\varphi_i)^2
=\eta\left(\varphi_i\cdot\delta_\mathbb{R}^G(\varphi_i)-\delta_\mathbb{R}^G(\wedge^2\varphi_i)\right).
\end{gather*}
In this example $KR$-theory can tell apart two dif\/ferent group actions, while $K$-theory cannot.
\end{Example}

\begin{Example}
Let $(G,\sigma_G)=(Sp(2m),\Id)$.
Then $R(Sp(2m))\cong\mathbb{Z}[\sigma_{2m}^1,\sigma_{2m}^2, \dots,\sigma_{2m}^m]$,
where $\sigma_{2m}^i$ is the class of the irreducible representation with highest weight $L_1+L_2+\dots+L_i$.
Note that $\sigma_{2m}^{2k-1}\in RH(Sp(2m),\Id , \mathbb{R})$, $\sigma_{2m}^{2k}\in RR(Sp(2m),\Id , \mathbb{R})$.
Moreover, $\sigma_{2m}^i+\wedge^{i-2}
\sigma_{2m}
=\wedge^i
\sigma_{2m}$ for $1\leq i\leq m$.
The equivariant $KR$-theory $KR^*_{(Sp(2m),\Id)}(Sp(2m),\Id)$ is isomorphic to, as $KR^*_{(Sp(2m),
\Id)}(\pt)$-modules, the exterior algebra over $KR^*_{(Sp(2m),\Id)}(\pt)$ generated by
$\delta_\mathbb{H}^G(\sigma_{2m}^{2k-1})$ and $\delta_\mathbb{R}^G(\sigma_{2m}^{2k})$ for $1\leq k\leq m$
by Theorem~\ref{mainthm2}.
The restriction map
\begin{gather*}
i^*: \ KR_{(U(2m),
\sigma_\mathbb{H})}^*(U(2m),
\sigma_\mathbb{H})\to KR_{(Sp(2m),\Id)}^*(Sp(2m),\Id)
\end{gather*}
sends $\delta_{\mathbb{R}}^G(\wedge^{2k}
\sigma_{2m}
-\wedge^{2k-2}
\sigma_{2m})$ to $\delta_{\mathbb{R}}^G(\sigma_{2m}^{2k})$ and $\delta_{\mathbb{H}}^G(\wedge^{2k+1}
\sigma_{2m}
-\wedge^{2k-1}
\sigma_{2m})$ to $\delta_{\mathbb{H}}^G(\sigma_{2m}^{2k+1})$.
Applying~$i^*$ to the relevant equations in Corollary~\ref{krringunitary}, we get
\begin{gather*}
\delta_\mathbb{R}^G\big(\sigma_{2m}^{2k}\big)^2=\eta
\sum\limits_{i=1}^{2k}\left(\sigma_{2m}^{4k-2i}\cdot\delta^G_\mathbb{R}\big(\sigma_{2m}^{2i}\big)+
\sigma_{2m}^{4k-2i+1}\cdot\delta^G_\mathbb{H}\big(\sigma_{2m}^{2i-1}\big)\right),
\\[-1pt]
\delta_\mathbb{H}^G\big(\sigma_{2m}^{2k-1}\big)^2=\eta
\sum\limits_{i=1}^{2k-1}\left(\sigma_{2m}^{4k-2-2i}\cdot\delta^G_\mathbb{R}\big(\sigma_{2m}^{2i}\big)+
\sigma_{2m}^{4k-1-2i}\cdot\delta^G_\mathbb{H}\big(\sigma_{2m}^{2i-1}\big)\right).
\end{gather*}
\end{Example}

\begin{Example}
Let $(G,\sigma_G)=(G_2,\Id)$.
Then $RR(G_2)\cong\mathbb{Z}[\sigma_1, \sigma_2]$, where $\sigma_1$ and $\sigma_2$
are the classes of irreducible representations of dimensions~7 and~14, respectively.
Note that both~$\sigma_1$ and~$\sigma_2$ are in $RR(G_2,\Id, \mathbb{R})$, and that $\wedge^2
\sigma_1=\sigma_1+\sigma_2$, $\wedge^2\sigma_2=\sigma_1^3-\sigma_1^2-2\sigma_1\sigma_2-\sigma_1$.
The equivariant $KR$-theory $KR^*_{(G_2,\Id
)}(G_2,\Id)$ is isomorphic to, as $KR^*_{(G_2,\Id)}(\pt)$-modules, the exterior algebra over
$KR^*_{(G_2,\Id)}(\pt)$ generated by $\delta_\mathbb{R}^G(\sigma_1)$ and $\delta_\mathbb{R}^G(\sigma_2)$,
by Theorem~\ref{equivstrthm}.
Using Theorem~\ref{equivkrthysquare} and Proposition~\ref{eqderivative}, we have
\begin{gather*}
\delta_\mathbb{R}^G(\sigma_1)^2=\eta\big((\sigma_1-1)\cdot\delta_\mathbb{R}^G(\sigma_1)+\delta_\mathbb{R}^G(\sigma_2)\big),
\qquad
\delta_\mathbb{R}^G(\sigma_2)^2
=\eta\big((\sigma_1^2-1)\cdot\delta_\mathbb{R}^G(\sigma_1)+\sigma_2\cdot\delta_\mathbb{R}^G(\sigma_2)\big).
\end{gather*}
\end{Example}

\subsection*{Acknowledgements}

The author would like to thank Professor Reyer Sjamaar for suggesting this problem, pains\-takingly proofreading the
manuscript, his patient guidance and encouragement throughout the course of this project.
He also thanks the referees for their critical comments, pointing out the relevance of the work~\cite{Bo} and a~mistake
in the def\/inition of $\varphi(d\rho)$ in~\cite{BZ} to him.

\pdfbookmark[1]{References}{ref}
\LastPageEnding

\end{document}